\documentclass[11pt,a4paper]{article}%
\usepackage{amscd}
\usepackage{amsmath}
\usepackage{graphicx}
\usepackage{amsfonts}
\usepackage{amssymb}
\usepackage{xcolor}
\usepackage{amsbsy,amsthm}
\usepackage{pifont}
\usepackage{fancyhdr,ifthen}
\usepackage[colorlinks=true, bookmarksnumbered=true, bookmarksopen=true,
bookmarksopenlevel=3, pdfstartview=FitH, linkcolor=cyan, pdfmenubar=true,
pdftoolbar=true, bookmarks=true,citecolor=cyan, urlcolor=magenta,
filecolor=cyan,menucolor=black,plainpages=false,pdfpagelabels]{hyperref}%
\DeclareMathOperator{\dv}{div}

\newtheorem{theorem}{Theorem}[section]
\newtheorem{lemma}[theorem]{Lemma}

\newtheorem{proposition}[theorem]{Proposition}
\newtheorem{definition}[theorem]{Definition}
\newtheorem{remark}[theorem]{Remark}

\numberwithin{equation}{section}

\begin{document}

\title{{\huge {On the well-posedness of the Hall-MHD system in}}\\{\huge {a critical setting of Besov-Morrey type}}\vspace{0.7cm} }
\author{{\Large {Lucas C. F. Ferreira\thanks{State University of Campinas (Unicamp),
IMECC-Department of Mathematics, Rua S\'{e}rgio Buarque de Holanda, CEP
13083-859, Campinas, SP, Brazil. Email:\ lcff@ime.unicamp.br (corresponding
author).} \ \ and \ Rafael P. da Silva\thanks{Federal University of Technology
- Paran\'{a} (UTFPR), Department of Mathematics, CEP 86300-000, Corn\'{e}lio
Proc\'{o}pio-PR, Brazil. Email:\ rpsilva@utfpr.edu.br.}}}\vspace{0.5cm}}
\date{}
\maketitle

\begin{abstract}
In this paper, we address the 3D incompressible Hall-magnetohydrodynamic
system (Hall-MHD). Our objective is to provide local and global well-posedness
results for initial velocity $u_{0}$, magnetic field $B_{0}$ and the current
$J_{0}:=\nabla\times B_{0}$ in a new critical framework, namely critical
Besov-Morrey spaces. These spaces combine typical characteristics from both
Besov and Morrey spaces allowing a broader framework that encompasses the
regularity properties inherent in Besov spaces with the Morrey space
structure. Compared to previous works in Sobolev and Besov spaces, our
approach accommodates a broader class of initial data, ensuring the
construction of a unique solution over time.



\

\end{abstract}

\pagestyle{fancy} \fancyhf{} \renewcommand{\headrulewidth}{0pt}
\chead{\ifthenelse{\isodd{\value{page}}}{Lucas C. F. Ferreira and Rafael P. da Silva}{Hall-MHD system in Besov-Morrey spaces}}
\rhead{\thepage}

\section{Introduction}

We consider the 3D Hall-magnetohydrodynamics system (Hall-MHD)%
\begin{equation}
\left\{
\begin{array}
[c]{rclll}%
\partial_{t}u+(u\cdot\nabla)u+\nabla\phi-(\nabla\times B)\times B & = &
\mu\Delta u, & \text{in} & \,\,(x,t)\in\mathbb{R}^{3}\times(0,\infty);\\
\partial_{t}B-\nabla\times((u-h(\nabla\times B))\times B) & = & \nu\Delta B, &
\text{in} & \,\,(x,t)\in\mathbb{R}^{3}\times(0,\infty);\\
\dv u=\dv B & = & 0, & \text{in} & \,\,(x,t)\in\mathbb{R}^{3}\times
(0,\infty);\\
u(x,0)=u_{0}(x)\text{ and }B(x,0) & = & B_{0}(x), & \text{in} & x\in
\mathbb{R}^{3},
\end{array}
\right.  \label{MHD}%
\end{equation}
where the unknown vector fields $u=u(x,t)$ and $B=B(x,t)$ represent the
velocity field and the magnetic field, respectively, and the scalar function
$\phi=\phi(x,t)$ represent the scalar pressure. The constants $\mu,\nu>0$
stand for the fluid viscosity and the magnetic resistivity, respectively. In
turn, the dimensionless parameter $h>0$ provides a measure of the intensity of
the Hall effect in the fluid. Moreover, we assume that the initial conditions
$u_{0}$ and $B_{0}$ are divergence-free, that is, $\dv u_{0}=\dv B_{0}=0$.

System (\ref{MHD}) emerges in modeling electrically conducting fluids
presenting diverse applications in physics such as geo-dynamo, neutron stars
and magnetic reconnection in plasmas (see \cite{Acher-Liu}, \cite{Lighthill}
and references therein). Comparing to the standard MHD system (see
\cite{Sermange}), we have the Hall electric field $hJ\times B$ with the
current $J=\nabla\times B,$ which in turn generates the Hall term
$h\nabla\times((\nabla\times B)\times B)$ in the left hand side of
(\ref{MHD})$_{2}$. This term entails an inherent difficulty in directly
applying the arguments developed for the study of the Navier-Stokes equations
and the classical MHD equations, and other related fluid dynamics models,
owing to their specific mathematical structures and constraints. For more
details on the physical background of (\ref{MHD}), we refer the reader to
\cite{Balbus-Terquem, Forbes, Huba, Mininni, Shalybkov-Urpin, Wardle}.

Compared to other fluid dynamics models, the Hall-MHD system has been
relatively underexplored in mathematical subjects such as well-posedness,
regularity, and stability of solutions, despite its significance in physics.
In what follows, we review some works from the literature that are related to
our purposes. In their pioneering work \cite{Acher-Liu}, Acheritogaray
\textit{et al.} provided a formal derivation of the Hall-MHD system from a
two-fluids Euler-Maxwell system for electrons and ions. Alternatively, they
put forth a kinetic formulation with a quasineutral coupling between a
Fokker-Planck equation for the ions and\ a suitable set of fluid equations for
the electrons. Moreover, they obtained global-in-time existence of weak
solutions in $L^{2}([0,1]^{3})$ by means of Galerkin method. From that point
on, Chae, Degond and Liu \cite{Chae-Degond-Liu} considered initial data
$u_{0}$ and $B_{0}$ in Sobolev spaces $H^{s}$ and established global existence
of weak solutions as well as the local well-posedness with $s>5/2$. Further
research in the context of weak solutions was conducted by Dumas and Sueur
\cite{Dumas-Sueur}. The global existence of strong solutions under smallness
condition in the initial data, as well as Serrin type continuation criteria
for smooth solutions, were obtained by Chae and Lee \cite{Chae-Lee} and Ye
\cite{Ye}. In \cite{Ahmad-Zhou}, Ahmad \textit{et al. }noted that the
potential blow-up of smooth solutions can be managed entirely in terms of the
velocity. Additionally, results on well-posedness of strong solutions for less
regular initial data in Sobolev or Besov spaces were established in
\cite{Benvenutti-Ferreira}, \cite{Wu-Yu-Tang} and \cite{Wan-Zhou}, and a
fundamental issue, which consists in the convergence as $h\rightarrow0$ to the
MHD system with no Hall-term, was discussed in \cite{Wan-Zhou-2}. In
\cite{Li-Yu-Zhu}, Li, Yu and Zhu showed examples of smooth data with
arbitrarily large norms in $L^{\infty}$ which generate unique global
solutions. More recently, in \cite{Danchin-Tan} and \cite{Danchin-Tan-2},
Danchin and Tan established local and global well-posedness results for
initial data in critical Sobolev and Besov spaces, respectively. In
\cite{Liu-Tan}, Liu and Tan extended some results of \cite{Danchin-Tan-2} to a
larger class of critical Besov spaces. Furthermore, in a very recent paper,
Fujii \cite{Fujii} obtained an existence-uniqueness result for global mild
solutions with the initial data in the space $\dot{B}_{p,\infty}^{\frac{3}%
{p}-1}\times(\dot{B}_{p,\infty}^{\frac{3}{p}-1}\cap L^{\infty})$,
$p\in(3,\infty)$, and assuming a suitable smallness condition. Alongside other
pertinent works on the subject, Nakasato \cite{Nakasato} established global
well-posedness for the Hall-MHD system around a constant equilibrium state
within the framework of critical Fourier-Besov spaces. Still studying
solutions as a perturbation from a constant equilibrium state, Kawashima,
Nakasato and Ogawa \cite{Kawashima-Nakasato-Ogawa} considered the compressible
case and proved the global well-posedness for initial data in critical $L^{2}%
$-Besov spaces. Zhang \cite{Zhang} achieved results on the local and global
well-posedness of strong solutions for initial velocity $u_{0}\in H^{\frac
{1}{2}+\sigma}$ with $0<\sigma<2$ and initial magnetic field $B_{0}\in
H^{\frac{3}{2}}$. An, Chen and Han \cite{An-Chen-Han} considered the problem
with variable density in critical Besov spaces and Han, Hu and Lai
\cite{Han-Hu-Lai} established global well-posedness for the compressible
Hall-MHD system with initial data in Besov spaces.

In this paper, we consider a new functional setting for the analysis of the
well-posedness of the Hall-MHD system, namely the homogeneous Besov-Morrey
spaces (abbreviated as BM-spaces). These spaces are denoted by $\mathcal{N}%
_{p,q,r}^{s}$ and belong to the category of Besov spaces with the underlying
space being the so-called Morrey spaces $\mathcal{M}_{q}^{p}$. We develop a
global-in-time existence and uniqueness theory for (\ref{MHD}) with initial
data $u_{0},B_{0}$ and $J_{0}:=\nabla\times B_{0}$ (the current) belonging to
$\mathcal{N}_{p,q,1}^{\frac{3}{p}-1}$ ($1\leq q\leq p<\infty$), under suitable
size conditions on the weak norm of the initial data. Comparing with the
standard homogeneous Besov space $\dot{B}_{p,r}^{s}$, we have that
$\mathcal{N}_{p,q,r}^{s}$ is strictly larger than it when $q<p$, $r\in
\lbrack1,\infty]$ and $s\in\mathbb{R},$ which allows us to consider a new
class of initial data for the well-posedness of (\ref{MHD}) (see Remark
\ref{paper-relevance} for more details). Moreover, if only $J_{0}$ is
sufficiently small, we are able to show a local version of the well-posedness
result. Our findings are built upon a suitable formulation of the Hall-MHD
system proposed by Danchin and Tan \cite{Danchin-Tan-2}, which uses the
current $J$ and bears similarities to the incompressible Navier-Stokes equations.

It is worth mentioning that Besov-Morrey spaces were introduced by Kozono and
Yamazaki \cite{Kozono-Yamazaki} in the 90s as a suitable extension that
combines typical characteristics from both Besov and Morrey spaces. This
merging was pivotal as it allowed for a broader framework that encompasses the
regularity properties inherent in Besov spaces with the Morrey space
structure. These spaces were initially employed in \cite{Kozono-Yamazaki} to
analyze the Navier-Stokes system by means of the so-called Kato approach (see
also \cite{Mazzucato}) and have been particularly influential in the analysis
of other fluid dynamics models and PDEs; see, e.g., \cite{Bie-Wang-Yao,
Duarte-Ferreira-Villa, Ferreira-Perez, Ferreira-Postigo, Xu-Tan, Yang-Fu-Sun}
and references therein.

Unlike the classical MHD system ($h=0$), the system under study lacks scaling
invariance due to the simultaneous presence of the Hall term in (\ref{MHD}%
)$_{2}$ and the Lorentz force in (\ref{MHD})$_{1}$. This coexistence
necessitates a precise understanding of critical regularity. For that,
following \cite{Danchin-Tan-2} and decoupling the system by making $u\equiv0$
and $B\equiv0$ in (\ref{MHD}), we are lead to the Navier-Stokes and Hall
systems, which are given by%

\begin{equation}
\left\{
\begin{array}
[c]{l}%
\partial_{t}u-\mu\Delta u+(u\cdot\nabla)u+\nabla\phi=0,\\
\dv u=0,
\end{array}
\right.  \label{aux-NS100}%
\end{equation}
and
\begin{equation}
\partial_{t}B-\nu\Delta B+h\nabla\times((\nabla\times B)\times B)=0,
\label{aux-Hall100}%
\end{equation}
respectively. These two systems present the respective scaling maps
\[
(u(x,t),\phi(x,t))\rightarrow(\lambda u(\lambda x,\lambda^{2}t),\lambda
^{2}\phi(\lambda x,\lambda^{2}t))\text{ and }B(x,t)\rightarrow B(\lambda
x,\lambda^{2}t).
\]
Thus, for $h>0\,\ $it would be natural to think that $u$ and $\nabla B$ have
the same level of homogeneity and then the same regularity. In turn, for the
classical MHD system ($h=0$)$,$ the fields $u$ and $B$ have the same homogeneity.

To make things compatible in the scaling sense, Danchin and Tan
\cite{Danchin-Tan-2} considered the current function $J=\nabla\times B$. In
fact, using that $\dv B=0,$ one arrives at $\Delta B=-(\nabla\times J)$ and
then the field $B$ can be written as%

\begin{equation}
B=(-\Delta)^{-1}\nabla\times J=curl^{-1}J, \label{B-identity}%
\end{equation}
where $curl^{-1}$ stands for the multiplier operator with symbol $i|\xi
|^{-2}\xi\times.$ From this, one gets the following system in terms of the
unknowns $u,B,J$%

\begin{align}
\partial_{t}u+(u\cdot\nabla)u-\mu\Delta u+\nabla\phi &  =J\times
B,\label{Ext.HMHD1}\\
\dv u =\dv B  &  =0,\label{Ext.HMHD2}\\
\partial_{t}B-\nabla\times((u-hJ)\times B)-\nu\Delta B  &
=0,\label{Ext.HMHD3}\\
\partial_{t}J-\nabla\times(\nabla\times((u-hJ)\times curl^{-1}J))-\nu\Delta J
&  =0, \label{Ext.HMHD4}%
\end{align}
which is called extended Hall-MHD system. Note that system (\ref{Ext.HMHD1}%
)-(\ref{Ext.HMHD4}) presents the scaling map
\begin{equation}
(u(x,t),B(x,t),J(x,t))\rightarrow(u,B,J)_{\lambda}:=\lambda(u(\lambda
x,\lambda^{2}t),B(\lambda x,\lambda^{2}t),J(\lambda x,\lambda^{2}t)),\text{
for }\lambda>0, \label{scal-ExHMHD-1}%
\end{equation}
which has the advantage of having the same homogeneity as the incompressible
Navier-Stokes equations.

In what follows, we present a well-posedness analysis in the framework of
critical BM-spaces for (\ref{Ext.HMHD1})-(\ref{Ext.HMHD4}). In this sense,
before stating our main results, we need to introduce some functional spaces
that will be used constantly throughout our work. For $1\leq q\leq p<\infty$
and $T>0$, we define
\[
X_{p}(T)=\left\{  u\in C([0,T];\mathcal{N}_{p,q,1}^{\frac{3}{p}-1})\cap
L^{1}(0,T;\mathcal{N}_{p,q,1}^{\frac{3}{p}+1});\text{ }\dv u=0\text{ in
}\mathcal{S}^{^{\prime}}(\mathbb{R}^{3})\right\}  ,
\]
equipped with the norm
\[
\Vert u\Vert_{X_{p}(T)}=\Vert u\Vert_{L^{\infty}(0,T;\mathcal{N}%
_{p,q,1}^{\frac{3}{p}-1})}+\Vert u\Vert_{L^{1}(0,T;\mathcal{N}_{p,q,1}%
^{\frac{3}{p}+1})}.
\]
Also, consider the global version of the above space by
\[
X_{p}=\left\{  u\in BC([0,\infty);\mathcal{N}_{p,q,1}^{\frac{3}{p}-1})\cap
L^{1}(0,\infty;\mathcal{N}_{p,q,1}^{\frac{3}{p}+1});\text{ }\dv u=0\text{ in
}\mathcal{S}^{^{\prime}}(\mathbb{R}^{3})\right\}  ,
\]
equipped with the norm
\[
\Vert u\Vert_{X_{p}}=\Vert u\Vert_{L^{\infty}(0,\infty;\mathcal{N}%
_{p,q,1}^{\frac{3}{p}-1})}+\Vert u\Vert_{L^{1}(0,\infty;\mathcal{N}%
_{p,q,1}^{\frac{3}{p}+1})}.
\]

Now we are in position to state our main results.

\begin{theorem}
\label{Result1} Let $1\leq q\leq p<\infty$, $u_{0},B_{0}\in\mathcal{N}%
_{p,q,1}^{\frac{3}{p}-1}$ with $\dv u_{0}=\dv B_{0}=0,$ and $J_{0}%
=\nabla\times B_{0}\in\mathcal{N}_{p,q,1}^{\frac{3}{p}-1}$. Then,
$\exists\delta=\delta(p,q,\mu,\nu)>0$ and an existence time $T>0$ such that,
if
\[
h\Vert J_{0}\Vert_{\mathcal{N}_{p,q,1}^{\frac{3}{p}-1}}<\delta,
\]
then system (\ref{MHD}) admits a unique local-in-time solution $(u,B)\in
X_{p}(T)\times X_{p}(T)$ with $J=\nabla\times B\in X_{p}(T)$. Moreover, the
solution depends continuously on the initial data and the estimates%
\[
h\Vert J\Vert_{X_{p}(T)}\leq c\left\{  h\Vert J_{0}\Vert_{\mathcal{N}%
_{p,q,1}^{\frac{3}{p}-1}}\right\}
\]
and
\[
\Vert u\Vert_{X_{p}(T)}+\Vert B\Vert_{X_{p}(T)}+h\Vert J\Vert_{X_{p}(T)}\leq
c\left\{  \Vert u_{0}\Vert_{\mathcal{N}_{p,q,1}^{\frac{3}{p}-1}}+\Vert
B_{0}\Vert_{\mathcal{N}_{p,q,1}^{\frac{3}{p}-1}}+h\Vert J_{0}\Vert
_{\mathcal{N}_{p,q,1}^{\frac{3}{p}-1}}\right\}
\]
hold true for some constant $c=c(p,q,\mu,\nu)>0.$
\end{theorem}

\begin{theorem}
\label{Result2} Let $1\leq q\leq p<\infty$, $u_{0},B_{0}\in\mathcal{N}%
_{p,q,1}^{\frac{3}{p}-1}$ with $\dv u_{0}=\dv B_{0}=0$ and $J_{0}=\nabla\times
B_{0}\in\mathcal{N}_{p,q,1}^{\frac{3}{p}-1}$. Then, system (\ref{MHD}) admits
a unique global-in-time solution $(u,B)\in X_{p}\times X_{p}$ with
$J=\nabla\times B\in X_{p}$, provided that
\[
\Vert u_{0}\Vert_{\mathcal{N}_{p,q,1}^{\frac{3}{p}-1}}+\Vert B_{0}%
\Vert_{\mathcal{N}_{p,q,1}^{\frac{3}{p}-1}}+h\Vert J_{0}\Vert_{\mathcal{N}%
_{p,q,1}^{\frac{3}{p}-1}}<\delta,
\]
where $\delta=\delta(p,q,\mu,\nu)$ is a positive constant. Moreover, the
solution depends continuously on the initial data and we have the estimate
\[
\Vert u\Vert_{X_{p}}+\Vert B\Vert_{X_{p}}+h\Vert J\Vert_{X_{p}}\leq c\left\{
\Vert u_{0}\Vert_{\mathcal{N}_{p,q,1}^{\frac{3}{p}-1}}+\Vert B_{0}%
\Vert_{\mathcal{N}_{p,q,1}^{\frac{3}{p}-1}}+h\Vert J_{0}\Vert_{\mathcal{N}%
_{p,q,1}^{\frac{3}{p}-1}}\right\}  ,
\]
for some constant $c=c(\mu,\nu)>0$.
\end{theorem}

Some comments on Theorems \ref{Result1} and \ref{Result2} are in order.

\begin{remark}
\label{paper-relevance}

\begin{itemize}
\item[(i)] In \cite{Danchin-Tan-2}, the authors established, among other
results, the global well-posedness of Hall-MHD system with $(u_{0},B_{0}%
,J_{0})\in(\dot{B}_{p,1}^{\frac{3}{p}-1})^{3}$ under a smallness condition on
the corresponding vector norm, and the local-in-time well-posedness under a
smallness condition on the $J_{0}$ norm. As already pointed out more above,
the strict inclusion $\dot{B}_{p,r}^{s}\hookrightarrow\mathcal{N}_{p,q,r}^{s}$
holds for $q<p,$ $r\in\lbrack1,\infty]$, and $s\in\mathbb{R}$ (see \cite[p.
964]{Kozono-Yamazaki}). Then, Theorem \ref{Result1} and Theorem \ref{Result2}
extend those results of \cite{Danchin-Tan-2} for a larger class of initial
data. Moreover, we have that $\mathcal{N}_{p,q,1}^{\frac{3}{p}-1}%
\not \subset \dot{B}_{2,r}^{\frac{1}{2}}$ \ for $q<p$ and any $r\in
\lbrack1,\infty]$. Then, our results and those obtained in
\cite{Danchin-Tan-2} in the spaces $\dot{B}_{2,r}^{\frac{1}{2}}$ with
$r\in\lbrack1,\infty]$ provide different initial data classes for the
well-posedness of (\ref{MHD}).

\item[(ii)] In \cite{Fujii}, the author carried out an existence-uniqueness
analysis of mild solutions with initial data $(u_{0},B_{0})\in\dot
{B}_{p,\infty}^{\frac{3}{p}-1}\times(\dot{B}_{p,\infty}^{\frac{3}{p}-1}\cap
L^{\infty})$, having obtained a result in the range $3<p<\infty$, applicable
to sufficiently small initial data in the norm of space $\dot{B}%
_{\infty,\infty}^{-1}\times(\dot{B}_{\infty,\infty}^{-1}\cap\dot{B}%
_{\infty,\infty}^{0})$. Given that there is no inclusion relation between the
spaces $\dot{B}_{p,\infty}^{\frac{3}{p}-1}$ and $\mathcal{N}_{\tilde{p}%
,\tilde{q},1}^{\frac{3}{\tilde{p}}-1}$ with $\tilde{q}<\tilde{p}$, the
initial-data classes considered in Theorem \ref{Result2} and in \cite{Fujii}
are distinct.
\end{itemize}
\end{remark}

\begin{remark}
\label{Rem-Simp-1}Note that the rescaled triplet $(\tilde{u},\tilde{B}%
,\tilde{\phi})(x,t):=\frac{h}{\mu}\left(  u,B,\frac{h}{\mu}\phi\right)
\left(  hx,\frac{h^{2}}{\mu}t\right)  $ satisfies system (\ref{MHD}) with
coefficients $(1,1,\nu/\mu)$ provided that $(u,B,\phi)$ is also a solution of
the same system with coefficients $(\mu,\nu,h).$ Considering this fact and the
scaling invariance of the homogeneous Besov-Morrey norms, it suffices to prove
the statements assuming $\mu=1$ and $h=1$. Thus, for simplicity, we will adopt
this assumption in the proofs, see Sections \ref{Sec3} and \ref{Sec4} for more
details. For explanatory purposes, we will also suppose that the magnetic
resistivity $\nu=1$. The general case can be addressed by appropriately
adjusting the heat semi-group.
\end{remark}

Our paper is organized as follows. Section \ref{Sec2} is intended to provide
some useful notation, review definitions and properties about function spaces
such as Morrey and Besov-Morrey spaces that will be used throughout our work.
In Section \ref{Sec3}, we analyze local well-posedness and prove Theorem
\ref{Result1}. In Section \ref{Sec4}, we dedicate ourselves to the study of
global well-posedness and prove Theorem \ref{Result2}.

\section{Preliminaries}

\label{Sec2}This section focuses on compiling essential notations,
definitions, tools, and properties related to certain operators and functional
spaces that will be pertinent to our purposes. Throughout, the letter $C>0$
will denote a constant that may vary with each occurrence. We will commence
with a review of Morrey spaces, including their definition and key properties.
For additional details, refer to \cite{Kozono-Yamazaki, Rosenthal,
Sawano-Fazio-Hakim, Zhou}.

\begin{definition}
For $1\leq q\leq p<\infty$, the Morrey space $\mathcal{M}_{q}^{p}%
=\mathcal{M}_{q}^{p}(\mathbb{R}^{n})$ is the set of all function $u\in
L_{loc}^{q}(\mathbb{R}^{n})$ such that%

\[
\parallel u\parallel_{\mathcal{M}_{q}^{p}}:=\sup_{x_{0}\in\mathbb{R}^{n}%
,R>0}R^{n/p-n/q}\parallel u\parallel_{L^{q}(B_{R}(x_{0}))}<\infty,
\]
where $B_{R}(x_{0})$ stands for the closed ball in $\mathbb{R}^{n}$ with
radius $R>0$ and center $x_{0}$.
\end{definition}

\begin{remark}
The space $\mathcal{M}_{q}^{p}$ endowed with the norm $\Vert\ .\ \Vert
_{\mathcal{M}_{q}^{p}}$ is a Banach space. Furthermore, the continuous
inclusion $\mathcal{M}_{r}^{p}\subset\mathcal{M}_{q}^{p}$ holds for $1\leq
q<r\leq p<\infty$ and $\mathcal{M}_{p}^{p}=L^{p}(\mathbb{R}^{n})$. In the case
where $p=\infty,$ we can identify $\mathcal{M}_{q}^{p}$ with $L^{\infty
}(\mathbb{R}^{n})$, meaning $\mathcal{M}_{q}^{\infty}=L^{\infty}$.
\end{remark}

In the sequel, we revisit the Littlewood-Paley decomposition (see
\cite{Cannone-Meyer, Chemin, Lemarie} for further details). Specifically, let
$\psi,\varphi\in C^{\infty}(\mathbb{R}^{n})$ be functions such that
$supp\varphi\subset\left\{  \xi\in\mathbb{R}^{n}:3/4\leq|\xi|\leq8/3\right\}
$, $supp\psi\subset B_{4/3}(0)$, and satisfying%

\[
\displaystyle\sum_{j=-\infty}^{\infty}\varphi_{j}(\xi)=1,\forall\xi
\in\mathbb{R}^{n}\backslash\{0\}\quad\mbox{and}\quad\psi(\xi
)+\displaystyle\sum_{j=0}^{\infty}\varphi_{j}(\xi)=1,\forall\xi\in
\mathbb{R}^{n},
\]
where $\varphi_{j}(\xi):=\varphi(2^{-j}\xi),\forall j\in\mathbb{Z},$ and
$B_{4/3}(0)$ denotes the ball of radius $4/3$ centered at the origin. Under
these conditions, we have the Littlewood-Paley decomposition
\[
u=\displaystyle\sum_{j=-\infty}^{\infty}\mathcal{F}^{-1}[\varphi
_{j}\mathcal{F}[u]]=\displaystyle\sum_{j=-\infty}^{\infty}\mathcal{F}%
^{-1}[\varphi_{j}]\ast u,\text{ \ \ }\forall u\in\mathcal{S}_{h}^{^{\prime}},
\]
where $\mathcal{S}^{^{\prime}}=\mathcal{S}^{^{\prime}}(\mathbb{R}^{n})$ stands
for the space of tempered distributions and $\mathcal{S}_{h}^{^{\prime}%
}=\mathcal{S}_{h}^{^{\prime}}(\mathbb{R}^{n})$ is the space of all
$u\in\mathcal{S}^{^{\prime}}$ such that $\displaystyle\lim_{j\rightarrow
-\infty}\Vert\tilde{S}_{j}u\Vert_{L^{\infty}}=0$ with $\tilde{S}_{j}$ being
the low-frequency cut-off operator
\[
\tilde{S}_{j}u:=\mathcal{F}^{-1}[\psi(2^{-j}\xi)\mathcal{F}[u]]=\mathcal{F}%
^{-1}[\psi(2^{-j}\xi)]\ast u,\text{ for }j\in\mathbb{Z}.
\]

We are now ready to define homogeneous Besov-Morrey spaces and recall some of
their basic properties (see \cite{Bie-Wang-Yao, Chen-Chen, Kozono-Yamazaki,
Mazzucato}).

\begin{definition}
For $1\leq q\leq p<\infty,$ $r\in\lbrack1,\infty]$ and $s\in\mathbb{R}$ the
homogeneous Besov-Morrey space $\mathcal{N}_{p,q,r}^{s}=\mathcal{N}%
_{p,q,r}^{s}(\mathbb{R}^{n})$ is the set of all $u\in\mathcal{S}_{h}%
^{^{\prime}}$ such that
\[
\mathcal{F}^{-1}[\varphi_{j}]\ast u\in\mathcal{M}_{q}^{p},\forall
j\in\mathbb{Z},
\]
and
\[
\parallel u\parallel_{\mathcal{N}_{p,q,r}^{s}}:=\left\Vert \left\{
2^{sj}\parallel\mathcal{F}^{-1}[\varphi_{j}]\ast u\parallel_{\mathcal{M}%
_{q}^{p}}\right\}  _{j=-\infty}^{\infty}\right\Vert _{\ell^{r}}<\infty.
\]

\end{definition}

\begin{remark}

\begin{itemize}
\item[(i)] The space $\mathcal{N}_{p,q,r}^{s}$ is a Banach space and the
continuous inclusion $\mathcal{N}_{p,q,r}^{s}\subset\mathcal{S}_{h}^{^{\prime
}}$ holds. Moreover, from the inclusion relations in $\mathcal{M}_{q}^{p}$ and
$\ell^{r}(\mathbb{Z})$, it follows that $\mathcal{N}_{p,q_{1},r_{1}}%
^{s}\subset\mathcal{N}_{p,q_{2},r_{2}}^{s}$, for $1\leq q_{2}\leq q_{1}\leq
p<\infty$ and $1\leq r_{1}\leq r_{2}\leq\infty$.

\item[(ii)] For $1\leq q\leq p<\infty$ and $r,r_{1},r_{2}\in\lbrack1,\infty]$,
Besov-Morrey spaces satisfy the real interpolation relation
\[
\mathcal{N}_{p,q,r}^{s}=(\mathcal{N}_{p,q,r_{1}}^{s_{1}},\mathcal{N}%
_{p,q,r_{2}}^{s_{2}})_{\theta,r}%
\]
with $s_{1},s_{2}\in\mathbb{R}$ and $0<\theta<1$, such that $s_{1}\neq s_{2}$
and $s=(1-\theta)s_{1}+\theta s_{2}$.
\end{itemize}
\end{remark}

In the next two lemmas, we recall basic tools that will be useful in obtaining
a series of estimates in the Besov-Morrey setting. The first is a
Mikhlin-H\"{o}rmander type theorem in homogeneous Besov-Morrey spaces and
therefore it provides us with a instrument for treating a proper class of
multiplier operators between those spaces. The second one consists of a
product estimate, which is useful for handling the nonlinear terms we will be
working with. These results can be found in \cite{Kozono-Yamazaki} and
\cite{Mazzucato}, respectively.

\begin{lemma}
\label{Mikhlin-Hormander} Let $m,s\in\mathbb{R}$, $1\leq q\leq p<\infty,$
$r\in\lbrack1,\infty],$ and let $f\in C^{\lfloor n/2\rfloor+1}(\mathbb{R}%
^{n}\backslash\left\{  0\right\}  ).$ Suppose also that there exists $A>0$
such that we have the pointwise estimate
\[
|(\partial^{\alpha}/\partial\xi)f(\xi)|\leq A|\xi|^{m-|\alpha|},\text{
}\forall\xi\in\mathbb{R}^{n}\backslash\left\{  0\right\}  ,
\]
for all $\alpha\in(\mathbb{N}\cup\{0\})^{n}$ with $|\alpha|\leq\lfloor
n/2\rfloor+1$. Then, the multiplier operator $Pu=\mathcal{F}^{-1}%
[f(\xi)\mathcal{F}[u]]$ is continuous from $\mathcal{N}_{p,q,r}^{s}$ to
$\mathcal{N}_{p,q,r}^{s-m}$ and satisfies the estimate
\[
\Vert Pu\Vert_{\mathcal{N}_{p,q,r}^{s-m}}\leq CA\Vert u\Vert_{\mathcal{N}%
_{p,q,r}^{s}},
\]
for some constant $C=C(m,n)>0.$
\end{lemma}

\begin{lemma}
\label{Mazzucato} The homogeneous Besov-Morrey space $\mathcal{N}%
_{p,q,1}^{\frac{n}{p}}(\mathbb{R}^{n})$, where $1\leq q\leq p<\infty$, is a
Banach algebra. That is, there exists a constant $C=C(p,q,n)>0$ such that
\[
\Vert uv\Vert_{\mathcal{N}_{p,q,1}^{\frac{n}{p}}}\leq C\Vert u\Vert
_{\mathcal{N}_{p,q,1}^{\frac{n}{p}}}\Vert v\Vert_{\mathcal{N}_{p,q,1}%
^{\frac{n}{p}}},\text{ }%
\]
for all $u,v\in\mathcal{N}_{p,q,1}^{\frac{n}{p}}$.
\end{lemma}

In our analysis we will need some regularity estimates for the heat semigroup
$\{e^{t\Delta}\}_{t\geq0}$ in homogeneous Besov-Morrey spaces which can be
found in the recent papers \cite{Nogayama-Sawano-1, Nogayama-Sawano-2,
Nogayama-Sawano-3}.

\begin{lemma}
\label{heat-kernel} Let $u_{0}\in\mathcal{N}_{p,q,1}^{s}$ and
$y(t):=e^{t\Delta}u_{0}$. Then, $y\in BC([0,\infty);\mathcal{N}_{p,q,1}%
^{s})\cap L^{1}(0,\infty;\mathcal{N}_{p,q,1}^{s+2})$ and we have the estimate
\[
\Vert y\Vert_{L^{\infty}\left(  0,\infty;\mathcal{N}_{p,q,1}^{s}\right)
}+\Vert y\Vert_{L^{1}\left(  0,\infty;\mathcal{N}_{p,q,1}^{s+2}\right)  }\leq
C\Vert u_{0}\Vert_{\mathcal{N}_{p,q,1}^{s}},
\]
for some constant $C=C(n)>0$. In particular, for all $T>0$, we have that $y\in
C([0,T];\mathcal{N}_{p,q,1}^{s})\cap L^{1}(0,T;\mathcal{N}_{p,q,1}^{s+2})$
and
\[
\Vert y\Vert_{L^{\infty}\left(  0,T;\mathcal{N}_{p,q,1}^{s}\right)  }+\Vert
y\Vert_{L^{1}\left(  0,T;\mathcal{N}_{p,q,1}^{s+2}\right)  }\leq C\Vert
u_{0}\Vert_{\mathcal{N}_{p,q,1}^{s}}.
\]

\end{lemma}

\begin{lemma}
\label{inhomogeneous-term} Let $f\in L^{1}\left(  0,T;\mathcal{N}_{p,q,1}%
^{s}\right)  $, for some finite time $T>0$, and $z(t):=\displaystyle\int
_{0}^{t}e^{(t-\tau)\Delta}f(\tau)d\tau,$ for $0\leq t\leq T$. Then, $z\in
C([0,T];\mathcal{N}_{p,q,1}^{s})\cap L^{1}(0,T;\mathcal{N}_{p,q,1}^{s+2})$ and
we have the estimate
\begin{equation}
\Vert z\Vert_{L^{\infty}\left(  0,T;\mathcal{N}_{p,q,1}^{s}\right)  }+\Vert
z\Vert_{L^{1}\left(  0,T;\mathcal{N}_{p,q,1}^{s+2}\right)  }\leq C\Vert
f\Vert_{L^{1}\left(  0,T;\mathcal{N}_{p,q,1}^{s}\right)  },
\label{aux-heat-01}%
\end{equation}
for some constant $C=C(n)>0$. Moreover, if $f\in L^{1}(0,\infty;\mathcal{N}%
_{p,q,1}^{s})$ then $z\in BC([0,\infty);\mathcal{N}_{p,q,1}^{s})\cap
L^{1}(0,\infty;\mathcal{N}_{p,q,1}^{s+2})$ and estimate (\ref{aux-heat-01})
remains valid when $T=\infty$.
\end{lemma}

Building on the two preceding lemmas, one can directly deduce the following proposition.

\begin{proposition}
\label{heat-eq-regularity} Consider the nonhomogeneous linear heat equation
\begin{align}
\partial_{t}u-\Delta u  &  =f,\text{ in }\mathbb{R}^{n}\times(0,\infty
),\label{heat-equation1}\\
u(x,0)  &  =u_{0},\text{ in }\mathbb{R}^{n}. \label{heat-equation2}%
\end{align}
If $u_{0}\in\mathcal{N}_{p,q,1}^{s}$ and $f\in L^{1}\left(  0,T;\mathcal{N}%
_{p,q,1}^{s}\right)  $ for some finite time $T>0$, then problem
(\ref{heat-equation1})-(\ref{heat-equation2}) admits a unique solution $u\in
C([0,T];\mathcal{N}_{p,q,1}^{s})\cap L^{1}\left(  0,T;\mathcal{N}%
_{p,q,1}^{s+2}\right)  .$ Moreover, for some constant $C=C(n)>0,$ we have the
estimate
\[
\Vert u\Vert_{L^{\infty}(0,T;\mathcal{N}_{p,q,1}^{s})}+\Vert u\Vert
_{L^{1}(0,T;\mathcal{N}_{p,q,1}^{s+2})}\leq C\left\{  \Vert u_{0}%
\Vert_{\mathcal{N}_{p,q,1}^{s}}+\Vert f\Vert_{L^{1}(0,T;\mathcal{N}%
_{p,q,1}^{s})}\right\}  .
\]
In addition, if $f\in L^{1}(0,\infty;\mathcal{N}_{p,q,1}^{s})$ then $u\in
BC([0,\infty);\mathcal{N}_{p,q,1}^{s})\cap L^{1}(0,\infty;\mathcal{N}%
_{p,q,1}^{s+2})$ and the above estimate still holds true when $T=\infty$.
\end{proposition}

We intend to approach the extended Hall-MHD system (\ref{Ext.HMHD1}%
)-(\ref{Ext.HMHD4}) by means of the Contraction Mapping Principle in the
Besov-Morrey setting. In order to streamline our analysis and avoid extensive
fixed-point computations, we will leverage the following corollary of the
fixed-point theorem. For related lemmas and a detailed discussion, please
refer to \cite{Lemarie}, \cite{Ferreira-Postigo} and
\cite{Ferreira-Villamizar}. The case $M>0$ will be used in proving Theorem
\ref{Result1}, while the case $M=0$ will be utilized in the proof of Theorem
\ref{Result2}.

\begin{lemma}
\label{auxiliar1} Let $(X,\left\Vert \cdot\right\Vert _{X})$ stand for a
Banach space and consider the bilinear operator $\mathcal{B}:X\times
X\longrightarrow X$ satisfying $\left\Vert \mathcal{B}(x_{1},x_{2})\right\Vert
\leq K\left\Vert x_{1}\right\Vert \left\Vert x_{2}\right\Vert ,$ $\forall
x_{1},x_{2}\in X$, for some constant $K>0.$ Suppose also that $\mathcal{L}%
:X\longrightarrow X$ is a bounded linear operator with norm $0\leq M<1.$ Then,
for each fixed $y\in X$ with $\Vert y\Vert_{X}<\displaystyle\frac{(1-M)^{2}%
}{4K},$ the equation $x=y+\mathcal{L}(x)+\mathcal{B}(x,x)$ has a unique
solution such that $\Vert x\Vert_{X}\leq\frac{1-M}{2K}.$ Moreover, we have
that%
\[
\Vert x\Vert_{X}\leq\frac{M-M^{2}}{2K}+2\Vert y\Vert_{X}.
\]

\end{lemma}

In what follows, inspired by \cite{Danchin-Tan-2}, we rewrite the system in a
way amenable to solution through a contraction argument. For that, we employ
certain vector identities. For any pair $v$ and $w$ of divergence-free vector
fields in $\mathbb{R}^{3}$, we have%

\begin{equation}
(w\cdot\nabla)v=\dv(v\otimes w) \label{v-grad-w}%
\end{equation}
and%

\begin{equation}
(\nabla\times w)\times w=(w\cdot\nabla)w-\nabla\left(  \frac{|w|^{2}}%
{2}\right)  . \label{rotwxw}%
\end{equation}
Denoting $\Lambda:=\phi+\displaystyle\frac{|B|^{2}}{2}$ and using the
identities (\ref{v-grad-w}) and (\ref{rotwxw}), equation (\ref{Ext.HMHD1}) can
be reformulated as
\begin{equation}
\partial_{t}u-\mu\Delta u+\nabla\Lambda=\dv(B\otimes B)-\dv(u\otimes u).
\label{use-leray}%
\end{equation}
Next, the identity
\begin{equation}
\nabla\times(w\times v)=(v\cdot\nabla)w-(w\cdot\nabla)v, \label{gradx(wxv)-V1}%
\end{equation}
combined with (\ref{v-grad-w}), leads us to
\begin{equation}
\nabla\times(w\times v)=\dv(w\otimes v)-\dv(v\otimes w). \label{gradx(wxv)-V2}%
\end{equation}
Then, applying the Leray projector $\mathcal{P}$ to the equation
(\ref{use-leray}) and making use of (\ref{gradx(wxv)-V2}) in equations
(\ref{Ext.HMHD3}) and (\ref{Ext.HMHD4}), the extended Hall-MHD system
(\ref{Ext.HMHD1})-(\ref{Ext.HMHD4}) can be reformulated as
\begin{align}
\partial_{t}u-\mu\Delta u  &  =\Pi_{a}(B,B)-\Pi_{a}(u,u),\label{MHDE1}\\
\partial_{t}B-\nu\Delta B  &  =\Pi_{b}(B,hJ-u),\label{MHDE2}\\
\partial_{t}J-\nu\Delta J  &  =\nabla\times\Pi_{b}(curl^{-1}%
J,hJ-u),\label{MHDE3}\\
(u,B,J)|_{t=0}  &  =(u_{0},B_{0},J_{0}), \label{MHDE4}%
\end{align}
where
\[
\Pi_{a}(v,w)=\frac{1}{2}\mathcal{P}(\dv(v\otimes w)+\dv(w\otimes v))
\]
and
\[
\Pi_{b}(v,w)=\dv(v\otimes w)-\dv(w\otimes v).
\]
Denoting $\Theta:=(u,B,J)$, the system (\ref{MHDE1})-(\ref{MHDE4}) can be
rewritten more compactly as
\begin{align}
\partial_{t}\Theta-\Delta_{\mu,\nu}\Theta &  =\Pi(\Theta,\Theta
),\label{System1}\\
\Theta|_{t=0}  &  =\Theta_{0}, \label{System2}%
\end{align}
where $\Theta_{0}=(u_{0},B_{0},J_{0})$, $\Delta_{\mu,\nu}\Theta=%
\begin{pmatrix}
\mu\Delta u\\
\nu\Delta B\\
\nu\Delta J
\end{pmatrix}
$ and $\Pi:\mathbb{R}^{3}\times\mathbb{R}^{3}\longrightarrow\mathbb{R}^{9}$ is
defined as
\[
\Pi(\Phi,\Psi)=%
\begin{pmatrix}
\Pi_{a}(\Phi_{2},\Psi_{2})-\Pi_{a}(\Phi_{1},\Psi_{1})\\
\Pi_{b}(\Phi_{2},h\Psi_{3}-\Psi_{1})\\
\nabla\times\Pi_{b}(curl^{-1}\Phi_{3},h\Psi_{3}-\Psi_{1})
\end{pmatrix}
,
\]
in which $\Phi=(\Phi_{1},\Phi_{2},\Phi_{3})$ and $\Psi=(\Psi_{1},\Psi_{2}%
,\Psi_{3})$.\newline

Considering the extended system over the initial one offers an advantage due
to its semi-linearity, contrasting with the quasi-linearity of the Hall-MHD
system for $(u,B)$. The quadratic terms in (\ref{MHDE1}) and (\ref{MHDE2})
closely resemble those of the incompressible Navier-Stokes equation. However,
the presence of the Hall-MHD term in (\ref{MHDE3}) necessitates a departure
from the traditional theory of Navier-Stokes equations. This is because the
differentiation in this case occurs outside the first variable of $\Pi_{b}$,
rather than within it. For further details, the reader is referred to
\cite{Danchin-Tan-2}.

We conclude this section by introducing some notations which will be used henceforth.%

\[
L^{p}\left(  0,\infty;\mathcal{N}_{p,q,r}^{s}\right)  =L^{p}(\mathcal{N}%
_{p,q,r}^{s}),\quad L^{p}\left(  0,T;\mathcal{N}_{p,q,r}^{s}\right)
=L_{T}^{p}(\mathcal{N}_{p,q,r}^{s}),\text{ for }T>0,
\]

\[
BC\left(  [0,\infty);\mathcal{N}_{p,q,r}^{s}\right)  =BC(\mathcal{N}%
_{p,q,r}^{s}),\quad C\left(  [0,T];\mathcal{N}_{p,q,r}^{s}\right)
=C_{T}(\mathcal{N}_{p,q,r}^{s}),\text{ for }T>0.
\]

\section{Local-in-time well-posedness}

\label{Sec3}In this section we prove the local well-posedness for the Hall-MHD
system which has been stated in Theorem \ref{Result1}. For that, we start by
showing the local-in-time existence for the extended system (\ref{System1}%
)-(\ref{System2}). In view of Remark \ref{Rem-Simp-1}, without loss of
generality, we also can assume $\mu=1,$ $h=1$ and $\nu=1.$

Given that we have a smallness condition only on $J_{0}$, the plan is to
utilize Lemma \ref{auxiliar1} and control the existence time according to the
solution $\Theta^{L}=(u^{L},B^{L},J^{L})$ of the homogeneous heat equation%

\begin{align}
\partial_{t}\Theta^{L}-\Delta\Theta^{L}  &  =0,\label{HE1}\\
\Theta^{L}|_{t=0}  &  =\Theta_{0}, \label{HE2}%
\end{align}
where $\Theta_{0}=(u_{0},B_{0},J_{0}).$ By Proposition
\ref{heat-eq-regularity}, for every $T>0$, we have that if $\Theta_{0}%
\in\mathcal{N}_{p,q,1}^{s}$ then $\Theta^{L}\in L_{T}^{1}(\mathcal{N}%
_{p,q,1}^{s+2})\cap C_{T}(\mathcal{N}_{p,q,1}^{s})$ and
\begin{equation}
\Vert\Theta^{L}\Vert_{L_{T}^{\infty}(\mathcal{N}_{p,q,1}^{s})}+\Vert\Theta
^{L}\Vert_{L_{T}^{1}(\mathcal{N}_{p,q,1}^{s+2})}\leq C\Vert\Theta_{0}%
\Vert_{\mathcal{N}_{p,q,1}^{s}}. \label{ThetaL-estimate}%
\end{equation}
In particular, we have
\begin{equation}
\Vert J^{L}\Vert_{L_{T}^{\infty}(\mathcal{N}_{p,q,1}^{s})}+\Vert J^{L}%
\Vert_{L_{T}^{1}(\mathcal{N}_{p,q,1}^{s+2})}\leq C\Vert J_{0}\Vert
_{\mathcal{N}_{p,q,1}^{s}}. \label{JL-estimate}%
\end{equation}
Furthermore, for all $u\in L_{T}^{\infty}(\mathcal{N}_{p,q,1}^{s})\cap
L_{T}^{1}(\mathcal{N}_{p,q,1}^{s+2})$, an interpolation argument and
H\"{o}lder inequality in time yield%

\begin{align}
\Vert u\Vert_{L_{T}^{2}(\mathcal{N}_{p,q,1}^{s+1})}  &  =\left\{  \int_{0}%
^{T}\Vert u(t)\Vert_{\mathcal{N}_{p,q,1}^{s+1}}^{2}dt\right\}  ^{\frac{1}{2}%
}\nonumber\\
&  \leq C\left\{  \int_{0}^{T}\Vert u(t)\Vert_{\mathcal{N}_{p,q,1}^{s}}\Vert
u(t)\Vert_{\mathcal{N}_{p,q,1}^{s+2}}dt\right\}  ^{\frac{1}{2}}\nonumber\\
&  \leq C\Vert u\Vert_{L_{T}^{\infty}(\mathcal{N}_{p,q,1}^{s})}^{\frac{1}{2}%
}\Vert u\Vert_{L_{T}^{1}(\mathcal{N}_{p,q,1}^{s+2})}^{\frac{1}{2}}.
\label{L2T-estimate}%
\end{align}
Then, $\Theta^{L}\in L_{T}^{2}(\mathcal{N}_{p,q,1}^{s+1})$ and, by the
dominated convergence theorem,%

\begin{equation}
\label{lim}\lim_{T\to0}\Vert\Theta^{L}\Vert_{L^{1}_{T}(\mathcal{N}%
_{p,q,1}^{s+2})}=\lim_{T\to0}\Vert\Theta^{L}\Vert_{L^{2}_{T}(\mathcal{N}%
_{p,q,1}^{s+1})}=0.
\end{equation}

Next, note that $\Theta$ is solution of (\ref{System1})-(\ref{System2}) in
$[0,T]\times\mathbb{R}^{3}$ if, and only if, $\Theta=\Theta^{L}+\tilde{\Theta
}$, where
\[
\tilde{\Theta}(t)=\int_{0}^{t}e^{(t-\tau)\Delta}(\Pi(\tilde{\Theta},\Theta
^{L})+\Pi(\Theta^{L},\tilde{\Theta})+\Pi(\tilde{\Theta},\tilde{\Theta}%
)+\Pi(\Theta^{L},\Theta^{L}))(\tau)d\tau,\text{ }\forall t\in\lbrack0,T].
\]
Thus, defining the bilinear operator
\[
\mathcal{B}(\Phi,\Psi)(t)=\int_{0}^{t}e^{(t-\tau)\Delta}\Pi(\Phi,\Psi
)(\tau)d\tau,
\]
taking $\tilde{y}=\mathcal{B}(\Theta^{L},\Theta^{L})$ and considering the
linear operator $\mathcal{L}(\Phi)=\mathcal{B}(\Phi,\Theta^{L})+\mathcal{B}%
(\Theta^{L},\Phi)$, our problem can be reformulated as the equation%

\begin{equation}
\label{tildeU}\tilde{\Theta}(t)=\tilde{y}(t)+ \mathcal{L}(\tilde{\Theta
})(t)+\mathcal{B}(\tilde{\Theta}, \tilde{\Theta})(t).
\end{equation}

Consider now the space $X_{T}=X_{p}(T)\times X_{p}(T)\times X_{p}(T)$ endowed
with the norm
\[
\Vert\Phi\Vert_{X_{T}}:=\Vert\Phi\Vert_{L_{T}^{\infty}(\mathcal{N}%
_{p,q,1}^{\frac{3}{p}-1})}+\Vert\Phi\Vert_{L_{T}^{1}(\mathcal{N}%
_{p,q,1}^{\frac{3}{p}+1})}=\Vert\Phi_{1}\Vert_{X_{p}(T)}+\Vert\Phi_{2}%
\Vert_{X_{p}(T)}+\Vert\Phi_{3}\Vert_{X_{p}(T)}.
\]
With the aim of applying Lemma \ref{auxiliar1} to equation (\ref{tildeU}), we
are going to provide some estimates that will be important not only throughout
this section, but also in the analysis of the global well-posedness.

\begin{lemma}
\label{div(vxw)} Let $v,w\in\mathcal{N}_{p,q,1}^{\frac{3}{p}}$ with $\dv v=\dv
w=0$. Then, we have the estimate
\begin{equation}
\Vert\dv(v\otimes w)\Vert_{\mathcal{N}_{p,q,1}^{\frac{3}{p}-1}}\leq C\Vert
v\Vert_{\mathcal{N}_{p,q,1}^{\frac{3}{p}}}\Vert w\Vert_{\mathcal{N}%
_{p,q,1}^{\frac{3}{p}}}, \label{aux-product-1}%
\end{equation}
for some constant $C=C(p,q)>0$.
\end{lemma}

\noindent{\textbf{Proof:}} First, by Lemma \ref{Mikhlin-Hormander}, we have
\[
\Vert\dv(v\otimes w)\Vert_{\mathcal{N}_{p,q,1}^{\frac{3}{p}-1}}\leq C\Vert
v\otimes w\Vert_{\mathcal{N}_{p,q,1}^{\frac{3}{p}}}.
\]
Also, Lemma \ref{Mazzucato} yields%
\[
\Vert v\otimes w\Vert_{\mathcal{N}_{p,q,1}^{\frac{3}{p}}}\leq C\Vert
v\Vert_{\mathcal{N}_{p,q,1}^{\frac{3}{p}}}\Vert w\Vert_{\mathcal{N}%
_{p,q,1}^{\frac{3}{p}}}.
\]
Combining these last two estimates, we arrive at (\ref{aux-product-1}).

\begin{flushright}
\ding{110}
\end{flushright}

\begin{lemma}
\label{div(curlxw)} Let $v,w\in\mathcal{N}_{p,q,1}^{\frac{3}{p}-1}%
\cap\mathcal{N}_{p,q,1}^{\frac{3}{p}+1}$ with $\dv v=\dv w=0$. Then, we have
the estimate
\[
\Vert\dv((curl^{-1}v)\otimes w)\Vert_{\mathcal{N}_{p,q,1}^{\frac{3}{p}}}\leq
C\Vert v\Vert_{\mathcal{N}_{p,q,1}^{\frac{3}{p}-1}}^{\frac{1}{2}}\Vert
w\Vert_{\mathcal{N}_{p,q,1}^{\frac{3}{p}-1}}^{\frac{1}{2}}\Vert v\Vert
_{\mathcal{N}_{p,q,1}^{\frac{3}{p}+1}}^{\frac{1}{2}}\Vert w\Vert
_{\mathcal{N}_{p,q,1}^{\frac{3}{p}+1}}^{\frac{1}{2}},
\]
for some constant $C=C(p,q)>0$.
\end{lemma}

\noindent{\textbf{Proof:}} Using vectorial identities, we can see that
\[
\dv((curl^{-1}v)\otimes w)=w\cdot\nabla(curl^{-1}v).
\]
Then, by Lemma \ref{Mazzucato} and Lemma \ref{Mikhlin-Hormander}, we obtain
that
\begin{align}
\Vert\dv((curl^{-1}v)\otimes w)\Vert_{\mathcal{N}_{p,q,1}^{\frac{3}{p}}}  &
=\Vert w\cdot\nabla(curl^{-1}v)\Vert_{\mathcal{N}_{p,q,1}^{\frac{3}{p}}%
}\nonumber\\
&  \leq C\Vert w\Vert_{\mathcal{N}_{p,q,1}^{\frac{3}{p}}}\Vert\nabla
(curl^{-1}v)\Vert_{\mathcal{N}_{p,q,1}^{\frac{3}{p}}}\nonumber\\
&  \leq C\Vert w\Vert_{\mathcal{N}_{p,q,1}^{\frac{3}{p}}}\Vert curl^{-1}%
v\Vert_{\mathcal{N}_{p,q,1}^{\frac{3}{p}+1}}\nonumber\\
&  \leq C\Vert v\Vert_{\mathcal{N}_{p,q,1}^{\frac{3}{p}}}\Vert w\Vert
_{\mathcal{N}_{p,q,1}^{\frac{3}{p}}}. \label{aux-interp-1}%
\end{align}

Now, since $\mathcal{N}_{p,q,1}^{\frac{3}{p}}=(\mathcal{N}_{p,q,1}^{\frac
{3}{p}-1},\mathcal{N}_{p,q,1}^{\frac{3}{p}+1})_{\frac{1}{2},1}$, an
interpolation argument leads us to
\[
\Vert v\Vert_{\mathcal{N}_{p,q,1}^{\frac{3}{p}}}\leq C\Vert v\Vert
_{\mathcal{N}_{p,q,1}^{\frac{3}{p}-1}}^{\frac{1}{2}}\Vert v\Vert
_{\mathcal{N}_{p,q,1}^{\frac{3}{p}+1}}^{\frac{1}{2}}\text{ \ and \ }\Vert
w\Vert_{\mathcal{N}_{p,q,1}^{\frac{3}{p}}}\leq C\Vert w\Vert_{\mathcal{N}%
_{p,q,1}^{\frac{3}{p}-1}}^{\frac{1}{2}}\Vert w\Vert_{\mathcal{N}%
_{p,q,1}^{\frac{3}{p}+1}}^{\frac{1}{2}},
\]
which together with (\ref{aux-interp-1}) gives the desired estimate.

\begin{flushright}
\ding{110}
\end{flushright}

\begin{lemma}
\label{div(vxcurl)} Let $v,w\in\mathcal{N}_{p,q,1}^{\frac{3}{p}-1}%
\cap\mathcal{N}_{p,q,1}^{\frac{3}{p}+1}$ with $\dv v=\dv w=0$. Then, we have
the estimate
\[
\Vert\dv(v\otimes(curl^{-1}w))\Vert_{\mathcal{N}_{p,q,1}^{\frac{3}{p}}}\leq
C\Vert v\Vert_{\mathcal{N}_{p,q,1}^{\frac{3}{p}+1}}\Vert w\Vert_{\mathcal{N}%
_{p,q,1}^{\frac{3}{p}-1}},
\]
for some constant $C=C(p,q)>0$.
\end{lemma}

\noindent{\textbf{Proof:}} Using $\dv(v\otimes(curl^{-1}w))=(curl^{-1}%
w)\cdot\nabla v$, Lemma \ref{Mazzucato} and Lemma \ref{Mikhlin-Hormander}, we
have that%

\begin{align*}
\Vert\dv(v\otimes(curl^{-1}w))\Vert_{\mathcal{N}_{p,q,1}^{\frac{3}{p}}}  &
=\Vert(curl^{-1}w)\cdot\nabla v\Vert_{\mathcal{N}_{p,q,1}^{\frac{3}{p}}}\\
&  \leq C\Vert curl^{-1}w\Vert_{\mathcal{N}_{p,q,1}^{\frac{3}{p}}}\Vert\nabla
v\Vert_{\mathcal{N}_{p,q,1}^{\frac{3}{p}}}\\
&  \leq C\Vert w\Vert_{\mathcal{N}_{p,q,1}^{\frac{3}{p}-1}}\Vert
v\Vert_{\mathcal{N}_{p,q,1}^{\frac{3}{p}+1}},
\end{align*}
as requested.

\begin{flushright}
\ding{110}
\end{flushright}

As a consequence of the above estimates, we have the following result.

\begin{lemma}
\label{local-bilinear-operator} Let $\Phi,\Psi\in X_{T}$ and $\mathcal{B}%
(\Phi,\Psi)(t)=\displaystyle\int_{0}^{t}e^{(t-\tau)\Delta}\Pi(\Phi,\Psi
)(\tau)d\tau$. Then, $\mathcal{B}$ is a continuous bilinear operator from
$X_{T}\times X_{T}$ to $X_{T}$.
\end{lemma}

\noindent{\textbf{Proof:}} Consider $v,w\in X_{p}(T)$. Then, making use of
Lemma \ref{Mikhlin-Hormander}, Lemma \ref{div(vxw)} and H\"{o}lder inequality,
we can estimate
\begin{align}
\Vert\Pi_{a}(v,w)\Vert_{L_{T}^{1}(\mathcal{N}_{p,q,1}^{\frac{3}{p}-1})}  &
=\int_{0}^{T}\Vert\Pi_{a}(v,w)(t)\Vert_{\mathcal{N}_{p,q,1}^{\frac{3}{p}-1}%
}dt\nonumber\\
&  \leq C\int_{0}^{T}\left\Vert (\dv(v\otimes w)+\dv(w\otimes
v))(t)\right\Vert _{\mathcal{N}_{p,q,1}^{\frac{3}{p}-1}}dt\nonumber\\
&  \leq C\int_{0}^{T}\left\Vert v(t)\right\Vert _{\mathcal{N}_{p,q,1}%
^{\frac{3}{p}}}\left\Vert w(t)\right\Vert _{\mathcal{N}_{p,q,1}^{\frac{3}{p}}%
}dt\nonumber\\
&  \leq C\left\{  \int_{0}^{T}\Vert v(t)\Vert_{\mathcal{N}_{p,q,1}^{\frac
{3}{p}}}^{2}dt\right\}  ^{\frac{1}{2}}\left\{  \int_{0}^{T}\Vert
w(t)\Vert_{\mathcal{N}_{p,q,1}^{\frac{3}{p}}}^{2}dt\right\}  ^{\frac{1}{2}%
}\nonumber\\
&  =C\Vert v\Vert_{L_{T}^{2}(\mathcal{N}_{p,q,1}^{\frac{3}{p}})}\Vert
w\Vert_{L_{T}^{2}(\mathcal{N}_{p,q,1}^{\frac{3}{p}})}. \label{estloc.1}%
\end{align}

Now, since $\Pi_{b}(v,w)=\dv(v\otimes w)-\dv(w\otimes v)$, a parallel argument
to the previous one leads us to the estimate%

\begin{equation}
\Vert\Pi_{b}(v,w)\Vert_{L_{T}^{1}(\mathcal{N}_{p,q,1}^{\frac{3}{p}-1})}\leq
C\Vert v\Vert_{L_{T}^{2}(\mathcal{N}_{p,q,1}^{\frac{3}{p}})}\Vert
w\Vert_{L_{T}^{2}(\mathcal{N}_{p,q,1}^{\frac{3}{p}})}. \label{estloc.2}%
\end{equation}
Furthermore, making use of Lemma \ref{Mikhlin-Hormander}, Lemma
\ref{div(curlxw)}, Lemma \ref{div(vxcurl)} and H\"{o}lder inequality, we
obtain that%

\begin{align}
&  \Vert\nabla\times\Pi_{b}(curl^{-1}v,w)\Vert_{L_{T}^{1}(\mathcal{N}%
_{p,q,1}^{\frac{3}{p}-1})}\nonumber\\
&  =\int_{0}^{T}\Vert(\nabla\times\Pi_{b}(curl^{-1}v,w))(t)\Vert
_{\mathcal{N}_{p,q,1}^{\frac{3}{p}-1}}dt\nonumber\\
&  \leq C\int_{0}^{T}\Vert\Pi_{b}(curl^{-1}v,w)(t)\Vert_{\mathcal{N}%
_{p,q,1}^{\frac{3}{p}}}dt\nonumber\\
&  =C\int_{0}^{T}\Vert(\dv((curl^{-1}v)\otimes w)-\dv(w\otimes(curl^{-1}%
v)))(t)\Vert_{\mathcal{N}_{p,q,1}^{\frac{3}{p}}}dt\nonumber\\
&  \leq C\left[  \int_{0}^{T}\Vert(\dv((curl^{-1}v)\otimes
w))(t)|_{\mathcal{N}_{p,q,1}^{\frac{3}{p}}}dt+\int_{0}^{T}\Vert(\dv(w\otimes
(curl^{-1}v)))(t)|_{\mathcal{N}_{p,q,1}^{\frac{3}{p}}}dt\right] \nonumber\\
&  \leq C\left[  \int_{0}^{T}\left\Vert v(t)\right\Vert _{\mathcal{N}%
_{p,q,1}^{\frac{3}{p}}}\left\Vert w(t)\right\Vert _{\mathcal{N}_{p,q,1}%
^{\frac{3}{p}}}dt+\int_{0}^{T}\left\Vert v(t)\right\Vert _{\mathcal{N}%
_{p,q,1}^{\frac{3}{p}-1}}\left\Vert w(t)\right\Vert _{\mathcal{N}%
_{p,q,1}^{\frac{3}{p}+1}}dt\right] \nonumber\\
&  \leq C\left[  \Vert v\Vert_{L_{T}^{2}(\mathcal{N}_{p,q,1}^{\frac{3}{p}}%
)}\Vert w\Vert_{L_{T}^{2}(\mathcal{N}_{p,q,1}^{\frac{3}{p}})}+\Vert
v\Vert_{L_{T}^{\infty}(\mathcal{N}_{p,q,1}^{\frac{3}{p}-1})}\Vert
w\Vert_{L_{T}^{1}(\mathcal{N}_{p,q,1}^{\frac{3}{p}+1})}\right]  .
\label{estloc.3}%
\end{align}

Consider now $\Phi=(\Phi_{1},\Phi_{2},\Phi_{3})$ and $\Psi=(\Psi_{1},\Psi
_{2},\Psi_{3})\in X_{T}$. By estimates (\ref{estloc.1}),(\ref{estloc.2}) and
(\ref{estloc.3}), we have that%

\begin{align}
&  \Vert\Pi(\Phi,\Psi)\Vert_{L_{T}^{1}(\mathcal{N}_{p,q,1}^{\frac{3}{p}-1}%
)}\nonumber\\
&  =\Vert\Pi_{a}(\Phi_{2},\Psi_{2})-\Pi_{a}(\Phi_{1},\Psi_{1})\Vert_{L_{T}%
^{1}(\mathcal{N}_{p,q,1}^{\frac{3}{p}-1})}+\Vert\Pi_{b}(\Phi_{2},\Psi_{3}%
-\Psi_{1})\Vert_{L_{T}^{1}(\mathcal{N}_{p,q,1}^{\frac{3}{p}-1})}\nonumber\\
&  +\Vert\nabla\times\Pi_{b}(curl^{-1}\Phi_{3},\Psi_{3}-\Psi_{1})\Vert
_{L_{T}^{1}(\mathcal{N}_{p,q,1}^{\frac{3}{p}-1})}\nonumber\\
&  \leq C\left\{  \Vert\Phi\Vert_{L_{T}^{2}(\mathcal{N}_{p,q,1}^{\frac{3}{p}%
})}\Vert\Psi\Vert_{L_{T}^{2}(\mathcal{N}_{p,q,1}^{\frac{3}{p}})}+\Vert\Phi
_{3}\Vert_{L_{T}^{\infty}(\mathcal{N}_{p,q,1}^{\frac{3}{p}-1})}\Vert\Psi
\Vert_{L_{T}^{1}(\mathcal{N}_{p,q,1}^{\frac{3}{p}+1})}\right\}  .
\label{Q-L2-estimate}%
\end{align}
Making use of estimate (\ref{L2T-estimate}) in (\ref{Q-L2-estimate}), we
arrive at%

\begin{align}
&  \Vert\Pi(\Phi,\Psi)\Vert_{L_{T}^{1}(\mathcal{N}_{p,q,1}^{\frac{3}{p}-1}%
)}\nonumber\\
&  \leq C\left\{  \Vert\Phi\Vert_{L_{T}^{\infty}(\mathcal{N}_{p,q,1}^{\frac
{3}{p}-1})}^{\frac{1}{2}}\Vert\Phi\Vert_{L_{T}^{1}(\mathcal{N}_{p,q,1}%
^{\frac{3}{p}+1})}^{\frac{1}{2}}\Vert\Psi\Vert_{L_{T}^{\infty}(\mathcal{N}%
_{p,q,1}^{\frac{3}{p}-1})}^{\frac{1}{2}}\Vert\Psi\Vert_{L_{T}^{1}%
(\mathcal{N}_{p,q,1}^{\frac{3}{p}+1})}^{\frac{1}{2}}+\Vert\Phi_{3}\Vert
_{L_{T}^{\infty}(\mathcal{N}_{p,q,1}^{\frac{3}{p}-1})}\Vert\Psi\Vert
_{L_{T}^{1}(\mathcal{N}_{p,q,1}^{\frac{3}{p}+1})}\right\} \nonumber\\
&  \leq C\Vert\Phi\Vert_{X_{T}}\Vert\Psi\Vert_{X_{T}}, \label{aux-bili-1001}%
\end{align}
which shows that, for every pair $(\Phi,\Psi)\in X_{T}\times X_{T}$, we have
$\Pi(\Phi,\Psi)\in L_{T}^{1}(\mathcal{N}_{p,q,1}^{\frac{3}{p}-1})$. Therefore,
it follows from Lemma \ref{inhomogeneous-term} that $\mathcal{B}(\Phi,\Psi)\in
X_{T}$ and satisfies the estimate
\begin{equation}
\Vert\mathcal{B}(\Phi,\Psi)\Vert_{X_{T}}=\Vert\mathcal{B}(\Phi,\Psi
)\Vert_{L_{T}^{\infty}(\mathcal{N}_{p,q,1}^{\frac{3}{p}-1})}+\Vert
\mathcal{B}(\Phi,\Psi)\Vert_{L_{T}^{1}(\mathcal{N}_{p,q,1}^{\frac{3}{p}+1}%
)}\leq C\Vert\Pi(\Phi,\Psi)\Vert_{L_{T}^{1}(\mathcal{N}_{p,q,1}^{\frac{3}%
{p}-1})}. \label{aux-bili-1002}%
\end{equation}
Inserting (\ref{aux-bili-1001}) into (\ref{aux-bili-1002}) yields%

\[
\Vert\mathcal{B}(\Phi,\Psi)\Vert_{X_{T}}\leq C\Vert\Phi\Vert_{X_{T}}\Vert
\Psi\Vert_{X_{T}},\forall\Phi,\Psi\in X_{T},
\]
and then we are done.

\begin{flushright}
\ding{110}
\end{flushright}

\begin{lemma}
\label{local-linear-operator} Let $\Phi\in X_{T}$ and $\mathcal{L}%
(\Phi)=\mathcal{B}(\Phi,\Theta^{L})+\mathcal{B}(\Theta^{L},\Phi)$. Then,
$\mathcal{L}$ is a continuous linear operator from $X_{T}$ to $X_{T}$ and
satisfies the inequality
\[
\left\Vert \mathcal{L}\right\Vert _{X_{T}\rightarrow X_{T}}\leq C\left\{
\Vert\Theta^{L}\Vert_{L_{T}^{2}(\mathcal{N}_{p,q,1}^{\frac{3}{p}})}%
+\Vert\Theta^{L}\Vert_{L_{T}^{1}(\mathcal{N}_{p,q,1}^{\frac{3}{p}+1})}+\Vert
J^{L}\Vert_{L_{T}^{\infty}(\mathcal{N}_{p,q,1}^{\frac{3}{p}-1})}\right\}  ,
\]
for some constant $C=C(p,q)>0$, where $\left\Vert \cdot\right\Vert
_{X_{T}\rightarrow X_{T}}$ stands for the operator norm in $X_{T}.$
\end{lemma}

\noindent{\textbf{Proof:}} In view of Lemma \ref{inhomogeneous-term} and
estimate (\ref{Q-L2-estimate}), we have that%

\begin{equation}
\Vert\mathcal{B}(\Phi,\Psi)\Vert_{X_{T}}\leq C\left\{  \Vert\Phi\Vert
_{L_{T}^{2}(\mathcal{N}_{p,q,1}^{\frac{3}{p}})}\Vert\Psi\Vert_{L_{T}%
^{2}(\mathcal{N}_{p,q,1}^{\frac{3}{p}})}+\Vert\Phi_{3}\Vert_{L_{T}^{\infty
}(\mathcal{N}_{p,q,1}^{\frac{3}{p}-1})}\Vert\Psi\Vert_{L_{T}^{1}%
(\mathcal{N}_{p,q,1}^{\frac{3}{p}+1})}\right\}  , \label{B-L2-estimate}%
\end{equation}
for all $\Phi,\Psi\in X_{T}$. Therefore, it is clear from estimate
(\ref{B-L2-estimate}) that%

\[
\Vert\mathcal{B}(\Phi,\Theta^{L})\Vert_{X_{T}}\leq C\left\{  \Vert\Phi
\Vert_{L_{T}^{2}(\mathcal{N}_{p,q,1}^{\frac{3}{p}})}\Vert\Theta^{L}%
\Vert_{L_{T}^{2}(\mathcal{N}_{p,q,1}^{\frac{3}{p}})}+\Vert\Phi_{3}\Vert
_{L_{T}^{\infty}(\mathcal{N}_{p,q,1}^{\frac{3}{p}-1})}\Vert\Theta^{L}%
\Vert_{L_{T}^{1}(\mathcal{N}_{p,q,1}^{\frac{3}{p}+1})}\right\}
\]
and
\[
\Vert\mathcal{B}(\Theta^{L},\Phi)\Vert_{X_{T}}\leq C\left\{  \Vert\Theta
^{L}\Vert_{L_{T}^{2}(\mathcal{N}_{p,q,1}^{\frac{3}{p}})}\Vert\Phi\Vert
_{L_{T}^{2}(\mathcal{N}_{p,q,1}^{\frac{3}{p}})}+\Vert J^{L}\Vert
_{L_{T}^{\infty}(\mathcal{N}_{p,q,1}^{\frac{3}{p}-1})}\Vert\Phi\Vert
_{L_{T}^{1}(\mathcal{N}_{p,q,1}^{\frac{3}{p}+1})}\right\}  .
\]
Thus, for all $\Phi\in X_{T}$, it follows that

{\fontsize{10}{10}\selectfont%
\begin{align*}
\Vert\mathcal{L}(\Phi)\Vert_{X_{T}}  &  \leq C\left\{  \Vert\Theta^{L}%
\Vert_{L_{T}^{2}(\mathcal{N}_{p,q,1}^{\frac{3}{p}})}\Vert\Phi\Vert_{L_{T}%
^{2}(\mathcal{N}_{p,q,1}^{\frac{3}{p}})}+\Vert\Theta^{L}\Vert_{L_{T}%
^{1}(\mathcal{N}_{p,q,1}^{\frac{3}{p}+1})}\Vert\Phi_{3}\Vert_{L_{T}^{\infty
}(\mathcal{N}_{p,q,1}^{\frac{3}{p}-1})}+\Vert J^{L}\Vert_{L_{T}^{\infty
}(\mathcal{N}_{p,q,1}^{\frac{3}{p}-1})}\Vert\Phi\Vert_{L_{T}^{1}%
(\mathcal{N}_{p,q,1}^{\frac{3}{p}+1})}\right\} \\
&  \leq C\left\{  \Vert\Theta^{L}\Vert_{L_{T}^{2}(\mathcal{N}_{p,q,1}%
^{\frac{3}{p}})}+\Vert\Theta^{L}\Vert_{L_{T}^{1}(\mathcal{N}_{p,q,1}^{\frac
{3}{p}+1})}+\Vert J^{L}\Vert_{L_{T}^{\infty}(\mathcal{N}_{p,q,1}^{\frac{3}%
{p}-1})}\right\}  \Vert\Phi\Vert_{X_{T}},
\end{align*}
}which concludes the proof.

\begin{flushright}
\ding{110}
\end{flushright}

With the assist of the above estimates, we now move on to the proof of Theorem
\ref{Result1}.

\noindent{\textbf{Proof of Theorem \ref{Result1}:}} Let $u_{0},B_{0}%
\in\mathcal{N}_{p,q,1}^{\frac{3}{p}-1}$ with $\dv u_{0}=\dv B_{0}=0$ and
$J_{0}=\nabla\times B_{0}\in\mathcal{N}_{p,q,1}^{\frac{3}{p}-1}$. Take
$\Theta_{0}=(u_{0},B_{0},J_{0})$ and consider the equation (\ref{tildeU}). By
Lemma \ref{local-bilinear-operator} we have that $\mathcal{B}:X_{T}\times
X_{T}\longrightarrow X_{T}$ is a continuous bilinear operator. Also, by Lemma
\ref{local-linear-operator}, $\mathcal{L}:X_{T}\longrightarrow X_{T}$ is a
continuous linear operator with the norm $M=\left\Vert \mathcal{L}\right\Vert
_{X_{T}\rightarrow X_{T}}$ satisfying%

\[
M\leq C\left\{  \Vert\Theta^{L}\Vert_{L_{T}^{2}(\mathcal{N}_{p,q,1}^{\frac
{3}{p}})}+\Vert\Theta^{L}\Vert_{L_{T}^{1}(\mathcal{N}_{p,q,1}^{\frac{3}{p}%
+1})}+\Vert J^{L}\Vert_{L_{T}^{\infty}(\mathcal{N}_{p,q,1}^{\frac{3}{p}-1}%
)}\right\}  ,
\]
for some constant $C>0$. For $T>0$ small enough, it is an immediate
consequence of the equality (\ref{lim}) and estimate (\ref{JL-estimate}) that%

\begin{equation}
M\leq C\Vert J_{0}\Vert_{\mathcal{N}_{p,q,1}^{\frac{3}{p}-1}}.
\label{M-estimate}%
\end{equation}
Then, for $\delta>0$ small enough, it follows that $M<1$ because $\Vert
J_{0}\Vert_{\mathcal{N}_{p,q,1}^{\frac{3}{p}-1}}<\delta$.

In addition, from estimate (\ref{B-L2-estimate}), we have that $\tilde
{y}=\mathcal{B}(\Theta^{L}, \Theta^{L})\in X_{T}$ and%

\[
\Vert\tilde{y}\Vert_{X_{T}}\leq C\left\{  \Vert\Theta^{L}\Vert_{L_{T}%
^{2}(\mathcal{N}_{p,q,1}^{\frac{3}{p}})}^{2}+\Vert J^{L}\Vert_{L_{T}^{\infty
}(\mathcal{N}_{p,q,1}^{\frac{3}{p}-1})}\Vert\Theta^{L}\Vert_{L_{T}%
^{1}(\mathcal{N}_{p,q,1}^{\frac{3}{p}+1})}\right\}  .
\]
Thus, for $T>0$ small enough, we obtain that%

\begin{equation}
\Vert\tilde{y}\Vert_{X_{T}}\leq C\Vert J_{0}\Vert_{\mathcal{N}_{p,q,1}%
^{\frac{3}{p}-1}}, \label{y-local-estimate}%
\end{equation}
and, taking $\delta>0$ small enough, we arrive at the estimate
\[
\Vert\tilde{y}\Vert_{X_{T}}<\frac{(1-M)^{2}}{4K},
\]
where $K=\left\Vert \mathcal{B}\right\Vert _{X_{T}\times X_{T}\rightarrow
X_{T}}$ is the norm of the bilinear operator $\mathcal{B}$ in $X_{T}.$

It follows from Lemma \ref{auxiliar1} that, for $T>0$ and $\delta>0$ small
enough, there exists a solution $\tilde{\Theta}\in X_{T}$ of equation
(\ref{tildeU}) that meets%

\[
\Vert\tilde{u}\Vert_{X_{p}(T)}+\Vert\tilde{B}\Vert_{X_{p}(T)}+\Vert\tilde
{J}\Vert_{X_{p}(T)}=\Vert\tilde{\Theta}\Vert_{X_{T}}\leq\frac{M-M^{2}}%
{2K}+2\Vert\tilde{y}\Vert_{X_{T}}\leq\frac{M}{2K}+2\Vert\tilde{y}\Vert_{X_{T}%
}\leq C\Vert J_{0}\Vert_{\mathcal{N}_{p,q,1}^{\frac{3}{p}-1}}.
\]
Consequently, there exists a solution $\Theta=\Theta^{L}+\tilde{\Theta}$ of
(\ref{System1})-(\ref{System2}) in $[0,T]\times\mathbb{R}^{3}$ and a constant
$c>0$ such that%

\[
\Vert\Theta\Vert_{X_{T}}=\Vert\Theta^{L} + \tilde{\Theta} \Vert_{X_{T}}%
\leq\Vert\Theta^{L} \Vert_{X_{T}}+\Vert\tilde{\Theta} \Vert_{X_{T}}\leq
c\Vert\Theta_{0}\Vert_{\mathcal{N}_{p,q,1}^{\frac{3}{p}-1}}
\]
and
\[
\Vert J \Vert_{X_{p}(T)}=\Vert J^{L} + \tilde{J} \Vert_{X_{p}(T)}\leq\Vert
J^{L} \Vert_{X_{p}(T)}+\Vert\tilde{J} \Vert_{X_{p}(T)}\leq c\Vert J_{0}%
\Vert_{\mathcal{N}_{p,q,1}^{\frac{3}{p}-1}}.
\]
This proves the local-in-time existence for the extended Hall-MHD system.

To complete the proof of local existence for the original Hall-MHD system, we
need to demonstrate that if $J_{0}=\nabla\times B_{0}$, then $J=\nabla\times
B$ in $(0,T]$ , thereby confirming that $(u,B)$ is indeed a local-in-time
solution of the system (\ref{MHD}) in the sense of distributions.

To do this, note that, since $curlB=\nabla\times B$, then $B=curl^{-1}%
(\nabla\times B)$. Replacing this equality on the equation (\ref{MHDE2}), we
have that%

\[
\partial_{t} B-\Delta B=\Pi_{b}(curl^{-1}(\nabla\times B),J-u).
\]
It follows from the above equality that%

\[
\partial_{t}(\nabla\times B)-\Delta(\nabla\times B)=\nabla\times\Pi
_{b}(curl^{-1}(\nabla\times B),J-u).
\]
Then, subtracting equation (\ref{MHDE3}) from the previous equality, we arrive at%

\begin{equation}
(\partial_{t}-\Delta)(\nabla\times B-J)=\nabla\times\Pi_{b}(curl^{-1}%
(\nabla\times B-J),J-u). \label{gradxB-J}%
\end{equation}
Consider now $0<T^{^{\prime}}\leq T$. Since $u,J\in L_{T}^{\infty}%
(\mathcal{N}_{p,q,1}^{\frac{3}{p}-1})\cap L_{T}^{1}(\mathcal{N}_{p,q,1}%
^{\frac{3}{p}+1})$, estimate (\ref{L2T-estimate}) allows us to state that
$u,J\in L_{T^{^{\prime}}}^{2}(\mathcal{N}_{p,q,1}^{\frac{3}{p}})$ and,
consequently, the same holds for $J-u$. Moreover, $u,J\in L_{T}^{\infty
}(\mathcal{N}_{p,q,1}^{\frac{3}{p}-1})$ implies that $u,J\in L_{T^{^{\prime}}%
}^{2}(\mathcal{N}_{p,q,1}^{\frac{3}{p}-1})$ and, since $B\in L_{T}^{\infty
}(\mathcal{N}_{p,q,1}^{\frac{3}{p}-1})\cap L_{T}^{1}(\mathcal{N}%
_{p,q,1}^{\frac{3}{p}+1})$, as well as $u$ and $J$, we have that $B\in
L_{T^{^{\prime}}}^{2}(\mathcal{N}_{p,q,1}^{\frac{3}{p}})$. So, by Lemma
\ref{Mikhlin-Hormander}, it follows that $\nabla\times B\in L_{T^{^{\prime}}%
}^{2}(\mathcal{N}_{p,q,1}^{\frac{3}{p}-1})$. In short, we have, for all
$0<T^{^{\prime}}\leq T$, that:

\begin{itemize}
\item {$J-u\in L_{T^{^{\prime}}}^{2}(\mathcal{N}_{p,q,1}^{\frac{3}{p}});$}

\item {$u,J\in L_{T^{^{\prime}}}^{2}(\mathcal{N}_{p,q,1}^{\frac{3}{p}-1});$}

\item {$\nabla\times B\in L_{T^{^{\prime}}}^{2}( \mathcal{N}_{p,q,1}^{\frac
{3}{p}-1}).$}
\end{itemize}

From definition of $\Pi_{b}$, we have%

\begin{align*}
\nabla\times\Pi_{b}(curl^{-1}(\nabla\times B-J),J-u)  &  =\nabla\times\left[
\dv(curl^{-1}(\nabla\times B-J)\otimes(J-u))\right. \\
&  -\left.  \dv((J-u)\otimes(curl^{-1}(\nabla\times B-J)))\right]  .
\end{align*}
Now, since $\nabla\times B-J\in L_{T^{^{\prime}}}^{2}(\mathcal{N}%
_{p,q,1}^{\frac{3}{p}-1})$, Lemma \ref{Mikhlin-Hormander} implies that
$curl^{-1}(\nabla\times B-J)\in L_{T^{^{\prime}}}^{2}(\mathcal{N}%
_{p,q,1}^{\frac{3}{p}})$. Also, using Lemma \ref{Mazzucato} and H\"{o}lder
inequality, we can estimate%
\begin{align*}
&  \Vert curl^{-1}(\nabla\times B-J)\otimes(J-u)\Vert_{L_{T^{^{\prime}}}%
^{1}(\mathcal{N}_{p,q,1}^{\frac{3}{p}})}\\
&  =\int_{0}^{T^{^{\prime}}}\Vert curl^{-1}(\nabla\times B-J)\otimes
(J-u)(t)\Vert_{\mathcal{N}_{p,q,1}^{\frac{3}{p}}}dt\\
&  \leq C\int_{0}^{T^{^{\prime}}}\Vert curl^{-1}(\nabla\times B-J)(t)\Vert
_{\mathcal{N}_{p,q,1}^{\frac{3}{p}}}\Vert(J-u)(t)\Vert_{\mathcal{N}%
_{p,q,1}^{\frac{3}{p}}}dt\\
&  \leq C\Vert curl^{-1}(\nabla\times B-J)\Vert_{L_{T^{^{\prime}}}%
^{2}(\mathcal{N}_{p,q,1}^{\frac{3}{p}})}\Vert J-u\Vert_{L_{T^{^{\prime}}}%
^{2}(\mathcal{N}_{p,q,1}^{\frac{3}{p}})}<\infty,
\end{align*}
which means that $curl^{-1}(\nabla\times B-J)\otimes(J-u)\in L_{T^{^{\prime}}%
}^{1}(\mathcal{N}_{p,q,1}^{\frac{3}{p}})$. This fact, together with Lemma
\ref{Mikhlin-Hormander}, implies that $\dv(curl^{-1}(\nabla\times
B-J)\otimes(J-u))\in L_{T^{^{\prime}}}^{1}(\mathcal{N}_{p,q,1}^{\frac{3}{p}%
-1})$. Also, with an analogous argument to the previous one, we obtain that
$\dv((J-u)\otimes(curl^{-1}(\nabla\times B-J)))\in L_{T^{^{\prime}}}%
^{1}(\mathcal{N}_{p,q,1}^{\frac{3}{p}-1})$. Consequently, applying again Lemma
\ref{Mikhlin-Hormander}, we arrive at $\nabla\times\Pi_{b}(curl^{-1}%
(\nabla\times B-J),J-u)\in L_{T^{^{\prime}}}^{1}(\mathcal{N}_{p,q,1}^{\frac
{3}{p}-2})$.

Note now that, by (\ref{gradxB-J}) and the hypothesis $J_{0}=\nabla\times B$,
we have that $\nabla\times B-J$ solves the problem%

\begin{align}
(\partial_{t}-\Delta)(\nabla\times B-J)  &  =\nabla\times\Pi_{b}%
(curl^{-1}(\nabla\times B-J),J-u),\label{problem:gradxB-J}\\
\nabla\times B-J|_{t=0}  &  =0.\nonumber
\end{align}
Therefore, it follows from Proposition \ref{heat-eq-regularity} that
$\nabla\times B-J\in C_{T^{^{\prime}}}(\mathcal{N}_{p,q,1}^{\frac{3}{p}%
-2})\cap L_{T^{^{\prime}}}^{1}(\mathcal{N}_{p,q,1}^{\frac{3}{p}})$ and%

\[
\Vert\nabla\times B-J\Vert_{L_{T^{^{\prime}}}^{\infty}(\mathcal{N}%
_{p,q,1}^{\frac{3}{p}-2})}+\Vert\nabla\times B-J\Vert_{L_{T^{^{\prime}}}%
^{1}(\mathcal{N}_{p,q,1}^{\frac{3}{p}})} \leq C \Vert\nabla\times\Pi
_{b}(curl^{-1}(\nabla\times B-J),J-u) \Vert_{L_{T^{^{\prime}}}^{1}%
(\mathcal{N}_{p,q,1}^{\frac{3}{p}-2})}.
\]
From this estimate, we get%

\begin{align}
&  \Vert(\nabla\times B-J)(T^{^{\prime}})\Vert_{\mathcal{N}_{p,q,1}^{\frac
{3}{p}-2}}+\int_{0}^{T^{^{\prime}}}\Vert(\nabla\times B-J)(t)\Vert
_{\mathcal{N}_{p,q,1}^{\frac{3}{p}}}dt\nonumber\\
&  \leq C\int_{0}^{T^{^{\prime}}}\Vert(\Pi_{b}(curl^{-1}(\nabla\times
B-J),J-u))(t)\Vert_{\mathcal{N}_{p,q,1}^{\frac{3}{p}-1}}dt\nonumber\\
&  \leq C\int_{0}^{T^{^{\prime}}}\Vert(curl^{-1}(\nabla\times B-J))(t)\Vert
_{\mathcal{N}_{p,q,1}^{\frac{3}{p}}}\Vert(J-u)(t)\Vert_{\mathcal{N}%
_{p,q,1}^{\frac{3}{p}}}dt\nonumber\\
&  \leq C\int_{0}^{T^{^{\prime}}}\Vert(\nabla\times B-J)(t)\Vert
_{\mathcal{N}_{p,q,1}^{\frac{3}{p}-1}}\Vert(J-u)(t)\Vert_{\mathcal{N}%
_{p,q,1}^{\frac{3}{p}}}dt. \label{Gronwall:gradxB-J}%
\end{align}
Since $\mathcal{N}_{p,q,1}^{\frac{3}{p}-1}=(\mathcal{N}_{p,q,1}^{\frac{3}%
{p}-2},\mathcal{N}_{p,q,1}^{\frac{3}{p}})_{\frac{1}{2},1}$, we can estimate%

\[
\Vert(\nabla\times B-J)(t)\Vert_{\mathcal{N}_{p,q,1}^{\frac{3}{p}-1}}\leq
C\Vert(\nabla\times B-J)(t)\Vert_{\mathcal{N}_{p,q,1}^{\frac{3}{p}-2}}%
^{\frac{1}{2}}\Vert(\nabla\times B-J)(t)\Vert_{\mathcal{N}_{p,q,1}^{\frac
{3}{p}}}^{\frac{1}{2}}.
\]
Then,

{\fontsize{10}{10}\selectfont%
\[
\Vert(\nabla\times B-J)(t) \Vert_{\mathcal{N}_{p,q,1}^{\frac{3}{p}-1}}%
\Vert(J-u)(t) \Vert_{\mathcal{N}_{p,q,1}^{\frac{3}{p}}}\leq C \Vert
(\nabla\times B-J)(t) \Vert^{\frac{1}{2}}_{\mathcal{N}_{p,q,1}^{\frac{3}{p}%
-2}}\Vert(\nabla\times B-J)(t) \Vert^{\frac{1}{2}}_{\mathcal{N}_{p,q,1}%
^{\frac{3}{p}}}\Vert(J-u)(t) \Vert_{\mathcal{N}_{p,q,1}^{\frac{3}{p}}}.
\]
} Applying Young inequality to the above estimate, we have, for all
$\epsilon>0$, that

{\fontsize{10}{10}\selectfont%
\[
\Vert(\nabla\times B-J)(t)\Vert_{\mathcal{N}_{p,q,1}^{\frac{3}{p}-1}}%
\Vert(J-u)(t)\Vert_{\mathcal{N}_{p,q,1}^{\frac{3}{p}}}\leq\epsilon\Vert
(\nabla\times B-J)(t)\Vert_{\mathcal{N}_{p,q,1}^{\frac{3}{p}}}+C\epsilon
^{-1}\Vert(\nabla\times B-J)(t)\Vert_{\mathcal{N}_{p,q,1}^{\frac{3}{p}-2}%
}\Vert(J-u)(t)\Vert_{\mathcal{N}_{p,q,1}^{\frac{3}{p}}}^{2}.
\]
} Replacing this estimate on (\ref{Gronwall:gradxB-J}), we arrive at%

\begin{align*}
&  \Vert(\nabla\times B-J)(T^{^{\prime}})\Vert_{\mathcal{N}_{p,q,1}^{\frac
{3}{p}-2}}+\int_{0}^{T^{^{\prime}}}\Vert(\nabla\times B-J)(t)\Vert
_{\mathcal{N}_{p,q,1}^{\frac{3}{p}}}dt\\
&  \leq C\left\{  \epsilon\int_{0}^{T^{^{\prime}}}\Vert(\nabla\times
B-J)(t)\Vert_{\mathcal{N}_{p,q,1}^{\frac{3}{p}}}dt+\epsilon^{-1}\int
_{0}^{T^{^{\prime}}}\Vert(J-u)(t)\Vert_{\mathcal{N}_{p,q,1}^{\frac{3}{p}}}%
^{2}\Vert(\nabla\times B-J)(t)\Vert_{\mathcal{N}_{p,q,1}^{\frac{3}{p}-2}%
}dt\right\}  .
\end{align*}
Thus, taking $\epsilon>0$ small enough, we obtain that%

\[
\Vert(\nabla\times B-J)(T^{^{\prime}})\Vert_{\mathcal{N}_{p,q,1}^{\frac{3}%
{p}-2}}+\frac{1}{2}\int_{0}^{T^{^{\prime}}}\Vert(\nabla\times B-J)(t)\Vert
_{\mathcal{N}_{p,q,1}^{\frac{3}{p}}}dt\leq C_{\epsilon}\int_{0}^{T^{^{\prime}%
}}\Vert(J-u)(t)\Vert_{\mathcal{N}_{p,q,1}^{\frac{3}{p}}}^{2}\Vert(\nabla\times
B-J)(t)\Vert_{\mathcal{N}_{p,q,1}^{\frac{3}{p}-2}}dt.
\]
Consequently,%

\[
\Vert(\nabla\times B-J)(T^{^{\prime}})\Vert_{\mathcal{N}_{p,q,1}^{\frac{3}%
{p}-2}}\leq C_{\epsilon}\int_{0}^{T^{^{\prime}}}\Vert(J-u)(t)\Vert
_{\mathcal{N}_{p,q,1}^{\frac{3}{p}}}^{2}\Vert(\nabla\times B-J)(t)\Vert
_{\mathcal{N}_{p,q,1}^{\frac{3}{p}-2}}dt.
\]
Then, making use of Gr\"{o}nwall inequality on the previous estimate, we
conclude that $\Vert(\nabla\times B-J)(T^{^{\prime}})\Vert_{\mathcal{N}%
_{p,q,1}^{\frac{3}{p}-2}}=0$. Since $0<T^{^{\prime}}\leq T$ is arbitrary, it
follows that $J=\nabla\times B$ in $(0,T]$.

Finally, we turn to the uniqueness property. For that, assume that
$(u^{1},B^{1})$ and $(u^{2},B^{2})$ are two solutions of the system
(\ref{MHD}) over the time interval $[0,T]$ with the same initial conditions.
Let $\Theta^{1}$ and $\Theta^{2}$ stand for the corresponding solutions to the
extended system (\ref{System1})-(\ref{System2}). Without losing generality, we
can take $\Theta^{2}$ as the solution constructed previously. Then, it follows
that%
\[
\Vert J^{2}\Vert_{X_{p}(T)}\leq c\Vert J_{0}\Vert_{\mathcal{N}_{p,q,1}%
^{\frac{3}{p}-1}},\text{ \ \ if }\Vert J_{0}\Vert_{\mathcal{N}_{p,q,1}%
^{\frac{3}{p}-1}}<\delta.
\]

Considering $\Phi=\Theta^{2}-\Theta^{1}$, we have that%

\begin{align*}
\partial_{t}\Phi-\Delta\Phi &  =\Pi(\Theta^{2},\Theta^{2})-\Pi(\Theta
^{1},\Theta^{1})\\
&  =\Pi(\Theta^{2},\Theta^{2})-\Pi(\Theta^{2},\Theta^{1})+\Pi(\Theta
^{2},\Theta^{1})-\Pi(\Theta^{1},\Theta^{1})\\
&  =\Pi(\Theta^{2},\Theta^{2}-\Theta^{1})+\Pi(\Theta^{2}-\Theta^{1},\Theta
^{1})\\
&  =\Pi(\Theta^{2},\Phi)+\Pi(\Phi,\Theta^{1}).
\end{align*}
Therefore, $\Phi$ is a solution of the problem%

\begin{align}
\partial_{t}\Phi-\Delta\Phi &  =\Pi(\Theta^{2},\Phi)+\Pi(\Phi,\Theta
^{1}),\label{UniSystem1}\\
\Phi|_{t=0}  &  =0, \label{UniSystem2}%
\end{align}
which yields%

\begin{equation}
\Phi(t)=\int_{0}^{t}e^{(t-\tau)\Delta}(\Pi(\Theta^{2},\Phi)+\Pi(\Phi
,\Theta^{1}))(\tau)d\tau=\mathcal{B}(\Theta^{2},\Phi)(t)+\mathcal{B}%
(\Phi,\Theta^{1})(t). \label{mildV}%
\end{equation}

Note now that, making use of Lemma \ref{inhomogeneous-term}, we have, for all
$T^{^{\prime}}\in(0,T]$, that%

\begin{align*}
&  \Vert\mathcal{B}(\Theta^{2},\Phi)\Vert_{L_{T^{^{\prime}}}^{\infty
}(\mathcal{N}_{p,q,1}^{\frac{3}{p}-1})}+\Vert\mathcal{B}(\Theta^{2},\Phi
)\Vert_{L_{T^{^{\prime}}}^{1}(\mathcal{N}_{p,q,1}^{\frac{3}{p}+1})}\\
&  \leq C\Vert\Pi(\Theta^{2},\Phi)\Vert_{L_{T^{^{\prime}}}^{1}(\mathcal{N}%
_{p,q,1}^{\frac{3}{p}-1})}\\
&  =C\left\{  \Vert\Pi_{a}(B^{2},\Phi_{2})-\Pi_{a}(u^{2},\Phi_{1}%
)\Vert_{L_{T^{^{\prime}}}^{1}(\mathcal{N}_{p,q,1}^{\frac{3}{p}-1})}+\Vert
\Pi_{b}(B^{2},\Phi_{3}-\Phi_{1})\Vert_{L_{T^{^{\prime}}}^{1}(\mathcal{N}%
_{p,q,1}^{\frac{3}{p}-1})}\right. \\
&  \left.  +\Vert\nabla\times\Pi_{b}(curl^{-1}J^{2},\Phi_{3}-\Phi_{1}%
)\Vert_{L_{T^{^{\prime}}}^{1}(\mathcal{N}_{p,q,1}^{\frac{3}{p}-1})}\right\}  .
\end{align*}
Then, by Lemma \ref{div(vxw)}, Lemma \ref{div(curlxw)} and Lemma
\ref{div(vxcurl)}, we get%

\begin{align}
&  \Vert\mathcal{B}(\Theta^{2},\Phi)\Vert_{L_{T^{^{\prime}}}^{\infty
}(\mathcal{N}_{p,q,1}^{\frac{3}{p}-1})}+\Vert\mathcal{B}(\Theta^{2},\Phi
)\Vert_{L_{T^{^{\prime}}}^{1}(\mathcal{N}_{p,q,1}^{\frac{3}{p}+1})}\nonumber\\
&  \leq C\left\{  \int_{0}^{T^{^{\prime}}}\Vert\Pi_{a}(B^{2},\Phi_{2}%
)(s)-\Pi_{a}(u^{2},\Phi_{1})(s)\Vert_{\mathcal{N}_{p,q,1}^{\frac{3}{p}-1}%
}+\int_{0}^{T^{^{\prime}}}\Vert\Pi_{b}(B^{2},\Phi_{3}-\Phi_{1})(s)\Vert
_{\mathcal{N}_{p,q,1}^{\frac{3}{p}-1}}ds\right. \nonumber\\
&  \left.  +\int_{0}^{T^{^{\prime}}}\Vert\Pi_{b}(curl^{-1}J^{2},\Phi_{3}%
-\Phi_{1})(s)\Vert_{\mathcal{N}_{p,q,1}^{\frac{3}{p}}}ds\right\} \nonumber\\
&  \leq C\left\{  \int_{0}^{T^{^{\prime}}}(\Vert B^{2}(s)\Vert_{\mathcal{N}%
_{p,q,1}^{\frac{3}{p}}}\Vert\Phi_{2}(s)\Vert_{\mathcal{N}_{p,q,1}^{\frac{3}%
{p}}}+\Vert u^{2}(s)\Vert_{\mathcal{N}_{p,q,1}^{\frac{3}{p}}}\Vert\Phi
_{1}(s)\Vert_{\mathcal{N}_{p,q,1}^{\frac{3}{p}}})ds\right. \nonumber\\
&  \left.  +\int_{0}^{T^{^{\prime}}}\Vert B^{2}(s)\Vert_{\mathcal{N}%
_{p,q,1}^{\frac{3}{p}}}\Vert(\Phi_{3}-\Phi_{1})(s)\Vert_{\mathcal{N}%
_{p,q,1}^{\frac{3}{p}}}ds+\int_{0}^{T^{^{\prime}}}\Vert\Pi_{b}(curl^{-1}%
J^{2},\Phi_{3}-\Phi_{1})(s)\Vert_{\mathcal{N}_{p,q,1}^{\frac{3}{p}}}ds\right\}
\nonumber\\
&  \leq C\left\{  \int_{0}^{T^{^{\prime}}}\Vert\Theta^{2}(s)\Vert
_{\mathcal{N}_{p,q,1}^{\frac{3}{p}}}\Vert\Phi(s)\Vert_{\mathcal{N}%
_{p,q,1}^{\frac{3}{p}}}ds+\int_{0}^{T^{^{\prime}}}\Vert\dv((curl^{-1}%
J^{2})\otimes(\Phi_{3}-\Phi_{1}))(s)\Vert_{\mathcal{N}_{p,q,1}^{\frac{3}{p}}%
}ds\right. \nonumber\\
&  \left.  +\int_{0}^{T^{^{\prime}}}\Vert\dv((\Phi_{3}-\Phi_{1})\otimes
(curl^{-1}J^{2}))(s)\Vert_{\mathcal{N}_{p,q,1}^{\frac{3}{p}}}ds\right\}
\nonumber\\
&  \leq C\left\{  \int_{0}^{T^{^{\prime}}}\Vert\Theta^{2}(s)\Vert
_{\mathcal{N}_{p,q,1}^{\frac{3}{p}}}\Vert\Phi(s)\Vert_{\mathcal{N}%
_{p,q,1}^{\frac{3}{p}}}ds+\int_{0}^{T^{^{\prime}}}\Vert J^{2}(s)\Vert
_{\mathcal{N}_{p,q,1}^{\frac{3}{p}}}\Vert(\Phi_{3}-\Phi_{1})(s)\Vert
_{\mathcal{N}_{p,q,1}^{\frac{3}{p}}}ds\right. \nonumber\\
&  \left.  +\int_{0}^{T^{^{\prime}}}\Vert J^{2}(s)\Vert_{\mathcal{N}%
_{p,q,1}^{\frac{3}{p}-1}}\Vert(\Phi_{3}-\Phi_{1})(s)\Vert_{\mathcal{N}%
_{p,q,1}^{\frac{3}{p}+1}}ds\right\} \nonumber\\
&  \leq C\left\{  \int_{0}^{T^{^{\prime}}}\Vert\Theta^{2}(s)\Vert
_{\mathcal{N}_{p,q,1}^{\frac{3}{p}}}\Vert\Phi(s)\Vert_{\mathcal{N}%
_{p,q,1}^{\frac{3}{p}}}ds+\int_{0}^{T^{^{\prime}}}\Vert J^{2}(s)\Vert
_{\mathcal{N}_{p,q,1}^{\frac{3}{p}-1}}\Vert\Phi(s)\Vert_{\mathcal{N}%
_{p,q,1}^{\frac{3}{p}+1}}ds\right\}  . \label{Uni-local-estimate1}%
\end{align}

Making use of interpolation and Young inequality in (\ref{Uni-local-estimate1}%
), we get, for all $\eta>0$, that%

\begin{align}
&  \Vert\mathcal{B}(\Theta^{2},\Phi)\Vert_{L_{T^{^{\prime}}}^{\infty
}(\mathcal{N}_{p,q,1}^{\frac{3}{p}-1})}+\Vert\mathcal{B}(\Theta^{2},\Phi
)\Vert_{L_{T^{^{\prime}}}^{1}(\mathcal{N}_{p,q,1}^{\frac{3}{p}+1})}\nonumber\\
&  \leq C\left\{  \int_{0}^{T^{^{\prime}}}\Vert\Theta^{2}(s)\Vert
_{\mathcal{N}_{p,q,1}^{\frac{3}{p}}}\Vert\Phi(s)\Vert_{\mathcal{N}%
_{p,q,1}^{\frac{3}{p}-1}}^{\frac{1}{2}}\Vert\Phi(s)\Vert_{\mathcal{N}%
_{p,q,1}^{\frac{3}{p}+1}}^{\frac{1}{2}}ds+\int_{0}^{T^{^{\prime}}}\Vert
J^{2}(s)\Vert_{\mathcal{N}_{p,q,1}^{\frac{3}{p}-1}}\Vert\Phi(s)\Vert
_{\mathcal{N}_{p,q,1}^{\frac{3}{p}+1}}ds\right\} \nonumber\\
&  \leq\eta\int_{0}^{T^{^{\prime}}}\Vert\Phi(s)\Vert_{\mathcal{N}%
_{p,q,1}^{\frac{3}{p}+1}}ds+C\eta^{-1}\int_{0}^{T^{^{\prime}}}\Vert\Theta
^{2}(s)\Vert_{\mathcal{N}_{p,q,1}^{\frac{3}{p}}}^{2}\Vert\Phi(s)\Vert
_{\mathcal{N}_{p,q,1}^{\frac{3}{p}-1}}ds+C\int_{0}^{T^{^{\prime}}}\Vert
J^{2}(s)\Vert_{\mathcal{N}_{p,q,1}^{\frac{3}{p}-1}}\Vert\Phi(s)\Vert
_{\mathcal{N}_{p,q,1}^{\frac{3}{p}+1}}ds\nonumber\\
&  \leq(\eta+C\Vert J^{2}\Vert_{L_{T^{^{\prime}}}^{\infty}(\mathcal{N}%
_{p,q,1}^{\frac{3}{p}-1})})\Vert\Phi\Vert_{L_{T^{^{\prime}}}^{1}%
(\mathcal{N}_{p,q,1}^{\frac{3}{p}+1})}+C\eta^{-1}\int_{0}^{T^{^{\prime}}}%
\Vert\Theta^{2}(s)\Vert_{\mathcal{N}_{p,q,1}^{\frac{3}{p}}}^{2}\Vert
\Phi(s)\Vert_{\mathcal{N}_{p,q,1}^{\frac{3}{p}-1}}ds.
\label{Uni-local-estimate2}%
\end{align}
Following the same steps taken to obtain (\ref{Uni-local-estimate1}), we get%

\begin{align*}
&  \Vert\mathcal{B}(\Phi,\Theta^{1})\Vert_{L_{T^{^{\prime}}}^{\infty
}(\mathcal{N}_{p,q,1}^{\frac{3}{p}-1})}+\Vert\mathcal{B}(\Phi,\Theta^{1}%
)\Vert_{L_{T^{^{\prime}}}^{1}(\mathcal{N}_{p,q,1}^{\frac{3}{p}+1})}\\
&  \leq C\left\{  \int_{0}^{T^{^{\prime}}}\Vert\Theta^{1}(s)\Vert
_{\mathcal{N}_{p,q,1}^{\frac{3}{p}}}\Vert\Phi(s)\Vert_{\mathcal{N}%
_{p,q,1}^{\frac{3}{p}}}ds+\int_{0}^{T^{^{\prime}}}\Vert\Theta^{1}%
(s)\Vert_{\mathcal{N}_{p,q,1}^{\frac{3}{p}+1}}\Vert\Phi(s)\Vert_{\mathcal{N}%
_{p,q,1}^{\frac{3}{p}-1}}ds\right\}  .
\end{align*}
Now, proceeding in the same way as we did to obtain (\ref{Uni-local-estimate2}%
), we have, for all $\eta>0$, that%

\begin{align}
&  \Vert\mathcal{B}(\Phi,\Theta^{1})\Vert_{L_{T^{^{\prime}}}^{\infty
}(\mathcal{N}_{p,q,1}^{\frac{3}{p}-1})}+\Vert\mathcal{B}(\Phi,\Theta^{1}%
)\Vert_{L_{T^{^{\prime}}}^{1}(\mathcal{N}_{p,q,1}^{\frac{3}{p}+1})}\nonumber\\
&  \leq\eta\Vert\Phi\Vert_{L_{T^{^{\prime}}}^{1}(\mathcal{N}_{p,q,1}^{\frac
{3}{p}+1})}+C\int_{0}^{T^{^{\prime}}}(\Vert\Theta^{1}(s)\Vert_{\mathcal{N}%
_{p,q,1}^{\frac{3}{p}+1}}+\eta^{-1}\Vert\Theta^{1}(s)\Vert_{\mathcal{N}%
_{p,q,1}^{\frac{3}{p}}}^{2})\Vert\Phi(s)\Vert_{\mathcal{N}_{p,q,1}^{\frac
{3}{p}-1}}ds. \label{Uni-local-estimate3}%
\end{align}
Then, it follows from equality (\ref{mildV}) and estimates
(\ref{Uni-local-estimate2}) and (\ref{Uni-local-estimate3}) that%

\begin{align*}
&  \Vert\Phi\Vert_{L_{T^{^{\prime}}}^{\infty}(\mathcal{N}_{p,q,1}^{\frac{3}%
{p}-1})}+\Vert\Phi\Vert_{L_{T^{^{\prime}}}^{1}(\mathcal{N}_{p,q,1}^{\frac
{3}{p}+1})}\\
&  \leq(2\eta+C\Vert J^{2}\Vert_{L_{T^{^{\prime}}}^{\infty}(\mathcal{N}%
_{p,q,1}^{\frac{3}{p}-1})})\Vert\Phi\Vert_{L_{T^{^{\prime}}}^{1}%
(\mathcal{N}_{p,q,1}^{\frac{3}{p}+1})}\\
&  +C\eta^{-1}\int_{0}^{T^{^{\prime}}}(\Vert\Theta^{1}(s)\Vert_{\mathcal{N}%
_{p,q,1}^{\frac{3}{p}+1}}+\Vert\Theta^{1}(s)\Vert_{\mathcal{N}_{p,q,1}%
^{\frac{3}{p}}}^{2}+\Vert\Theta^{2}(s)\Vert_{\mathcal{N}_{p,q,1}^{\frac{3}{p}%
}}^{2})\Vert\Phi(s)\Vert_{\mathcal{N}_{p,q,1}^{\frac{3}{p}-1}}ds.
\end{align*}
Since
\[
\Vert J^{2}\Vert_{L_{T^{^{\prime}}}^{\infty}(\mathcal{N}_{p,q,1}^{\frac{3}%
{p}-1})}\leq\Vert J^{2}\Vert_{L_{T}^{\infty}(\mathcal{N}_{p,q,1}^{\frac{3}%
{p}-1})}\leq\Vert J^{2}\Vert_{X_{p}(T)}\leq c\Vert J_{0}\Vert_{\mathcal{N}%
_{p,q,1}^{\frac{3}{p}-1}}\text{ \ \ and \ \ }\Vert J_{0}\Vert_{\mathcal{N}%
_{p,q,1}^{\frac{3}{p}-1}}<\delta,
\]
we have, for $\delta>0$ and $\eta>0$ small enough, that

{\fontsize{10}{10}\selectfont%
\[
\Vert\Phi\Vert_{L_{T^{^{\prime}}}^{\infty}(\mathcal{N}_{p,q,1}^{\frac{3}{p}%
-1})}+\frac{1}{2}\Vert\Phi\Vert_{L_{T^{^{\prime}}}^{1}(\mathcal{N}%
_{p,q,1}^{\frac{3}{p}+1})}\leq C\int_{0}^{T^{^{\prime}}} (\Vert\Theta
^{1}(s)\Vert_{\mathcal{N}_{p,q,1}^{\frac{3}{p}+1}}+\Vert\Theta^{1}(s)\Vert
^{2}_{\mathcal{N}_{p,q,1}^{\frac{3}{p}}}+\Vert\Theta^{2}(s)\Vert
^{2}_{\mathcal{N}_{p,q,1}^{\frac{3}{p}}})\Vert\Phi(s)\Vert_{\mathcal{N}%
_{p,q,1}^{\frac{3}{p}-1}}ds.
\]
} Consequently,

{\fontsize{10}{10}\selectfont%
\[
\Vert\Phi({T^{^{\prime}}})\Vert_{\mathcal{N}_{p,q,1}^{\frac{3}{p}-1}}\leq
\Vert\Phi\Vert_{L_{T^{^{\prime}}}^{\infty}(\mathcal{N}_{p,q,1}^{\frac{3}{p}%
-1})}\leq C\int_{0}^{T^{^{\prime}}}(\Vert\Theta^{1}(s)\Vert_{\mathcal{N}%
_{p,q,1}^{\frac{3}{p}+1}}+\Vert\Theta^{1}(s)\Vert_{\mathcal{N}_{p,q,1}%
^{\frac{3}{p}}}^{2}+\Vert\Theta^{2}(s)\Vert_{\mathcal{N}_{p,q,1}^{\frac{3}{p}%
}}^{2})\Vert\Phi(s)\Vert_{\mathcal{N}_{p,q,1}^{\frac{3}{p}-1}}ds.
\]
} Applying Gr\"{o}nwall inequality on the previous estimate, we get $\Vert
\Phi({T^{^{\prime}}})\Vert_{\mathcal{N}_{p,q,1}^{\frac{3}{p}-1}}=0,\forall
T^{^{\prime}}\in(0,T]$. Thus, $\Phi\equiv0$ and thereby $\Theta_{1}%
\equiv\Theta_{2}$, which concludes the proof.

\begin{flushright}
\ding{110}
\end{flushright}

\section{Global-in-time well-posedness}

\label{Sec4}This section aims to establish the proof of Theorem \ref{Result2}.
Following, in part, the reasoning used in the previous section, we shall first
prove the global existence for the extended system (\ref{System1}%
)-(\ref{System2}). The proof will be carried out with the help of Lemma
\ref{auxiliar1} with $M=0$. Moreover, taking into account Remark
\ref{Rem-Simp-1}, again we can consider $\mu=1,$ $h=1$ and $\nu=1 ,$ without
compromising generality.

In order to do this, consider the space $X=X_{p}\times X_{p}\times X_{p}$
endowed with the norm
\[
\Vert\Phi\Vert_{X}:=\Vert\Phi\Vert_{L^{\infty}(\mathcal{N}_{p,q,1}^{\frac
{3}{p}-1})}+\Vert\Phi\Vert_{L^{1}(\mathcal{N}_{p,q,1}^{\frac{3}{p}+1})}%
=\Vert\Phi_{1}\Vert_{X_{p}}+\Vert\Phi_{2}\Vert_{X_{p}}+\Vert\Phi_{3}%
\Vert_{X_{p}}.
\]
Take $y(t)=e^{t\Delta}\Theta_{0},$ $t\geq0,$ and the bilinear operator
$\mathcal{B}$ as in the previous section. By Duhamel principle, system
(\ref{System1})-(\ref{System2}) can be expressed in the following integral
form%
\begin{equation}
\Theta(t)=e^{t\Delta}\Theta_{0}+\int_{0}^{t}e^{(t-\tau)\Delta}\Pi
(\Theta,\Theta)(\tau)d\tau=y(t)+\mathcal{B}(\Theta,\Theta)(t).
\label{integral-equation}%
\end{equation}
In order to apply Lemma {\ref{auxiliar1}} to equation (\ref{integral-equation}%
), we need to prove the following properties for $\mathcal{B}(\cdot,\cdot).$

\begin{lemma}
\label{global-bilinear-operator} Let $\Phi,\Psi\in X$ and $\mathcal{B}%
(\Phi,\Psi)(t)=\displaystyle\int_{0}^{t}e^{(t-\tau)\Delta}\Pi(\Phi,\Psi
)(\tau)d\tau$. Then, $\mathcal{B}$ is a continuous bilinear operator from
$X\times X$ to $X$.
\end{lemma}

\noindent{\textbf{Proof:}} For $v,w\in X_{p}$, we can make use of Lemma
\ref{Mikhlin-Hormander}, Lemma \ref{div(vxw)} and afterwards interpolation and
H\"{o}lder inequality in the time variable in order to handle $\Pi_{a}(v,w)$
as follows:%

\begin{align}
\Vert\Pi_{a}(v,w)\Vert_{L^{1}(\mathcal{N}_{p,q,1}^{\frac{3}{p}-1})}  &
=\int_{0}^{\infty}\Vert\Pi_{a}(v,w)(t)\Vert_{\mathcal{N}_{p,q,1}^{\frac{3}%
{p}-1}}dt\nonumber\\
&  \leq C\int_{0}^{\infty}\left\Vert (\dv(v\otimes w)+\dv(w\otimes
v))(t)\right\Vert _{\mathcal{N}_{p,q,1}^{\frac{3}{p}-1}}dt\nonumber\\
&  \leq C\int_{0}^{\infty}\left\Vert v(t)\right\Vert _{\mathcal{N}%
_{p,q,1}^{\frac{3}{p}}}\left\Vert w(t)\right\Vert _{\mathcal{N}_{p,q,1}%
^{\frac{3}{p}}}dt\nonumber\\
&  \leq C\int_{0}^{\infty}\left\Vert v(t)\right\Vert _{\mathcal{N}%
_{p,q,1}^{\frac{3}{p}-1}}^{\frac{1}{2}}\left\Vert w(t)\right\Vert
_{\mathcal{N}_{p,q,1}^{\frac{3}{p}-1}}^{\frac{1}{2}}\left\Vert v(t)\right\Vert
_{\mathcal{N}_{p,q,1}^{\frac{3}{p}+1}}^{\frac{1}{2}}\left\Vert w(t)\right\Vert
_{\mathcal{N}_{p,q,1}^{\frac{3}{p}+1}}^{\frac{1}{2}}dt\nonumber\\
&  \leq C\Vert v\Vert_{L^{\infty}(\mathcal{N}_{p,q,1}^{\frac{3}{p}-1})}%
^{\frac{1}{2}}\Vert w\Vert_{L^{\infty}(\mathcal{N}_{p,q,1}^{\frac{3}{p}-1}%
)}^{\frac{1}{2}}\left\{  \int_{0}^{\infty}\Vert v(t)\Vert_{\mathcal{N}%
_{p,q,1}^{\frac{3}{p}+1}}dt\right\}  ^{\frac{1}{2}}\left\{  \int_{0}^{\infty
}\Vert w(t)\Vert_{\mathcal{N}_{p,q,1}^{\frac{3}{p}+1}}dt\right\}  ^{\frac
{1}{2}}\nonumber\\
&  \leq C\Vert v\Vert_{X_{p}}^{\frac{1}{2}}\Vert w\Vert_{X_{p}}^{\frac{1}{2}%
}\Vert v\Vert_{L^{1}(\mathcal{N}_{p,q,1}^{\frac{3}{p}+1})}^{\frac{1}{2}}\Vert
w\Vert_{L^{1}(\mathcal{N}_{p,q,1}^{\frac{3}{p}+1})}^{\frac{1}{2}}\nonumber\\
&  \leq C\Vert v\Vert_{X_{p}}\Vert w\Vert_{X_{p}}. \label{est.1}%
\end{align}

Now, since $\Pi_{b}(v,w)=\dv(v\otimes w)-\dv(w\otimes v)$, using the same
arguments as above, we obtain the estimate%

\begin{equation}
\Vert\Pi_{b}(v,w)\Vert_{L^{1}(\mathcal{N}_{p,q,1}^{\frac{3}{p}-1})}\leq C\Vert
v\Vert_{X_{p}}\Vert w\Vert_{X_{p}}. \label{est.2}%
\end{equation}
Furthermore, employing Lemma \ref{Mikhlin-Hormander}, we have that%

\begin{align}
&  \Vert\nabla\times\Pi_{b}(curl^{-1}v,w)\Vert_{L^{1}(\mathcal{N}%
_{p,q,1}^{\frac{3}{p}-1})}\nonumber\\
&  =\int_{0}^{\infty}\Vert(\nabla\times\Pi_{b}(curl^{-1}v,w))(t)\Vert
_{\mathcal{N}_{p,q,1}^{\frac{3}{p}-1}}dt\nonumber\\
&  \leq C\int_{0}^{\infty}\Vert\Pi_{b}(curl^{-1}v,w)(t)\Vert_{\mathcal{N}%
_{p,q,1}^{\frac{3}{p}}}dt\nonumber\\
&  =C\int_{0}^{\infty}\Vert(\dv((curl^{-1}v)\otimes w)-\dv(w\otimes
(curl^{-1}v)))(t)\Vert_{\mathcal{N}_{p,q,1}^{\frac{3}{p}}}dt\nonumber\\
&  \leq C\left[  \int_{0}^{\infty}(\dv((curl^{-1}v)\otimes w))(t)\Vert
_{\mathcal{N}_{p,q,1}^{\frac{3}{p}}}dt+\int_{0}^{\infty}\Vert(\dv(w\otimes
(curl^{-1}v)))(t)\Vert_{\mathcal{N}_{p,q,1}^{\frac{3}{p}}}dt\right]  .
\label{est.3-aux-1}%
\end{align}
In turn, using Lemma \ref{div(curlxw)}, Lemma \ref{div(vxcurl)} and H\"{o}lder
inequality in time, the R.H.S. of (\ref{est.3-aux-1}) can be estimated by%
\begin{align}
&  C\left[  \int_{0}^{\infty}\left\Vert v(t)\right\Vert _{\mathcal{N}%
_{p,q,1}^{\frac{3}{p}-1}}^{\frac{1}{2}}\left\Vert w(t)\right\Vert
_{\mathcal{N}_{p,q,1}^{\frac{3}{p}-1}}^{\frac{1}{2}}\left\Vert v(t)\right\Vert
_{\mathcal{N}_{p,q,1}^{\frac{3}{p}+1}}^{\frac{1}{2}}\left\Vert w(t)\right\Vert
_{\mathcal{N}_{p,q,1}^{\frac{3}{p}+1}}^{\frac{1}{2}}dt+\int_{0}^{\infty
}\left\Vert v(t)\right\Vert _{\mathcal{N}_{p,q,1}^{\frac{3}{p}-1}}\left\Vert
w(t)\right\Vert _{\mathcal{N}_{p,q,1}^{\frac{3}{p}+1}}dt\right] \nonumber\\
&  \leq C\left[  \Vert v\Vert_{X_{p}}^{\frac{1}{2}}\Vert w\Vert_{X_{p}}%
^{\frac{1}{2}}\int_{0}^{\infty}\left\Vert v(t)\right\Vert _{\mathcal{N}%
_{p,q,1}^{\frac{3}{p}+1}}^{\frac{1}{2}}\left\Vert w(t)\right\Vert
_{\mathcal{N}_{p,q,1}^{\frac{3}{p}+1}}^{\frac{1}{2}}dt+\Vert v\Vert_{X_{p}%
}\int_{0}^{\infty}\left\Vert w(t)\right\Vert _{\mathcal{N}_{p,q,1}^{\frac
{3}{p}+1}}dt\right] \nonumber\\
&  \leq C\left[  \Vert v\Vert_{X_{p}}^{\frac{1}{2}}\Vert w\Vert_{X_{p}}%
^{\frac{1}{2}}\Vert v\Vert_{L^{1}(\mathcal{N}_{p,q,1}^{\frac{3}{p}+1})}%
^{\frac{1}{2}}\Vert w\Vert_{L^{1}(\mathcal{N}_{p,q,1}^{\frac{3}{p}+1})}%
^{\frac{1}{2}}+\Vert v\Vert_{X_{p}}\Vert w\Vert_{L^{1}(\mathcal{N}%
_{p,q,1}^{\frac{3}{p}+1})}\right] \nonumber\\
&  \leq C\Vert v\Vert_{X_{p}}\Vert w\Vert_{X_{p}}. \label{est.3-aux-2}%
\end{align}
Consequently, combining estimate (\ref{est.3-aux-2}) and (\ref{est.3-aux-1})
enables us to arrive at%

\begin{equation}
\label{est.3}\Vert\nabla\times\Pi_{b}(curl^{-1}v,w)\Vert_{L^{1}(\mathcal{N}%
_{p,q,1}^{\frac{3}{p}-1})}\leq C \Vert v\Vert_{X_{p}}\Vert w\Vert_{X_{p}}.
\end{equation}

Consider now $\Phi=(\Phi_{1},\Phi_{2},\Phi_{3})$ and $\Psi=(\Psi_{1},\Psi
_{2},\Psi_{3})$ belonging to $X$. Then, by estimates (\ref{est.1}),
(\ref{est.2}) and (\ref{est.3}), we have that%

\begin{align*}
\Vert\Pi(\Phi,\Psi)\Vert_{L^{1}(\mathcal{N}_{p,q,1}^{\frac{3}{p}-1})}  &
=\Vert\Pi_{a}(\Phi_{2},\Psi_{2})-\Pi_{a}(\Phi_{1},\Psi_{1})\Vert
_{L^{1}(\mathcal{N}_{p,q,1}^{\frac{3}{p}-1})}\\
&  +\Vert\Pi_{b}(\Phi_{2},\Psi_{3}-\Psi_{1})\Vert_{L^{1}(\mathcal{N}%
_{p,q,1}^{\frac{3}{p}-1})}+\Vert\nabla\times\Pi_{b}(curl^{-1}\Phi_{3},\Psi
_{3}-\Psi_{1})\Vert_{L^{1}(\mathcal{N}_{p,q,1}^{\frac{3}{p}-1})}\\
&  \leq C\left\{  \Vert\Phi_{2}\Vert_{X_{p}}\Vert\Psi_{2}\Vert_{X_{p}}%
+\Vert\Phi_{1}\Vert_{X_{p}}\Vert\Psi_{1}\Vert_{X_{p}}+(\Vert\Phi_{2}%
\Vert_{X_{p}}+\Vert\Phi_{3}\Vert_{X_{p}})\Vert\Psi_{3}-\Psi_{1}\Vert_{X_{p}%
}\right\} \\
&  \leq C\Vert\Phi\Vert_{X}\Vert\Psi\Vert_{X},
\end{align*}
because $\Vert\Phi_{i}\Vert_{X_{p}}\leq\Vert\Phi\Vert_{X}\quad$and$\quad
\Vert\Psi_{i}\Vert_{X_{p}}\leq\Vert\Psi\Vert_{X},$ for $i=1,2,3.$ Thus, we get
the estimate%

\begin{equation}
\Vert\Pi(\Phi,\Psi)\Vert_{L^{1}(\mathcal{N}_{p,q,1}^{\frac{3}{p}-1})}\leq
C\Vert\Phi\Vert_{X}\Vert\Psi\Vert_{X}, \label{Q-estimate}%
\end{equation}
which shows that $\Pi(\Phi,\Psi)\in L^{1}(\mathcal{N}_{p,q,1}^{\frac{3}{p}%
-1})$ provided that $(\Phi,\Psi)\in X\times X.$ Therefore, it follows from
Lemma \ref{inhomogeneous-term} that $\mathcal{B}(\Phi,\Psi)\in X$ and
satisfies the estimate
\begin{equation}
\Vert\mathcal{B}(\Phi,\Psi)\Vert_{X}=\Vert\mathcal{B}(\Phi,\Psi)\Vert
_{L^{\infty}(\mathcal{N}_{p,q,1}^{\frac{3}{p}-1})}+\Vert\mathcal{B}(\Phi
,\Psi)\Vert_{L^{1}(\mathcal{N}_{p,q,1}^{\frac{3}{p}+1})}\leq C\Vert\Pi
(\Phi,\Psi)\Vert_{L^{1}(\mathcal{N}_{p,q,1}^{\frac{3}{p}-1})}.
\label{Q-estimate2}%
\end{equation}
Merging this last estimate with (\ref{Q-estimate}) promptly establishes the
desired continuity of $\mathcal{B}$.

\begin{flushright}
\ding{110}
\end{flushright}

We are now prepared to present the proof of Theorem \ref{Result2}, detailed as
follows.\vspace{0cm}

\noindent{\textbf{Proof of Theorem \ref{Result2}:}} Let $u_{0},B_{0}%
\in\mathcal{N}_{p,q,1}^{\frac{3}{p}-1}$ with $\dv u_{0}=\dv B_{0}=0$ and
$J_{0}=\nabla\times B_{0}\in\mathcal{N}_{p,q,1}^{\frac{3}{p}-1}$. Take
$\Theta_{0}=(u_{0},B_{0},J_{0})$ and consider the equation
(\ref{integral-equation}). From Lemma \ref{global-bilinear-operator} we have
that $\mathcal{B}:X\times X\longrightarrow X$ is a continuous bilinear
operator. Then, it follows immediately from Lemma \ref{heat-kernel} that $y\in
X$ and
\[
\Vert y\Vert_{X}\leq C\Vert\Theta_{0}\Vert_{\mathcal{N}_{p,q,1}^{\frac{3}%
{p}-1}}.
\]
Now, if%

\[
\Vert u_{0}\Vert_{\mathcal{N}_{p,q,1}^{\frac{3}{p}-1}}+\Vert B_{0}%
\Vert_{\mathcal{N}_{p,q,1}^{\frac{3}{p}-1}}+\Vert J_{0}\Vert_{\mathcal{N}%
_{p,q,1}^{\frac{3}{p}-1}}=\Vert\Theta_{0}\Vert_{\mathcal{N}_{p,q,1}^{\frac
{3}{p}-1}}<\frac{1}{4CK},
\]
where $K$ is the norm of operator $\mathcal{B}$, we have that $\Vert
y\Vert_{X}<\displaystyle\frac{1}{4K}$.

It follows from Lemma \ref{auxiliar1} with $M=0$ that there exists a solution
$\Theta\in X$ of equation (\ref{integral-equation}) that meets%

\[
\Vert u \Vert_{X_{p}}+\Vert B \Vert_{X_{p}}+\Vert J \Vert_{X_{p}}=\Vert
\Theta\Vert_{X}\leq2\Vert y \Vert_{X}\leq2C\Vert\Theta_{0}\Vert_{\mathcal{N}%
_{p,q,1}^{\frac{3}{p}-1}}.
\]
Therefore, there is a constant $c>0$ such that%

\[
\Vert u \Vert_{X_{p}}+\Vert B \Vert_{X_{p}}+\Vert J \Vert_{X_{p}}\leq
c\left\{  \Vert u_{0}\Vert_{\mathcal{N}_{p,q,1}^{\frac{3}{p}-1}}+\Vert
B_{0}\Vert_{\mathcal{N}_{p,q,1}^{\frac{3}{p}-1}}+\Vert J_{0}\Vert
_{\mathcal{N}_{p,q,1}^{\frac{3}{p}-1}}\right\}  ,
\]
which proves the global existence for the extended Hall-MHD system.

Note now that, since $u,B,J\in L^{\infty}(\mathcal{N}_{p,q,1}^{\frac{3}{p}%
-1})\cap L^{1}(\mathcal{N}_{p,q,1}^{\frac{3}{p}+1})$, an argument analogous to
the one made in estimate (\ref{L2T-estimate}) leads us to conclude that
$u,B,J\in L^{2}(\mathcal{N}_{p,q,1}^{\frac{3}{p}})$. Consequently, $J-u\in
L^{2}(\mathcal{N}_{p,q,1}^{\frac{3}{p}})$ and Lemma \ref{Mikhlin-Hormander}
implies that $\nabla\times B\in L^{2}(\mathcal{N}_{p,q,1}^{\frac{3}{p}-1})$.
In addition, $u,J\in L_{T}^{2}(\mathcal{N}_{p,q,1}^{\frac{3}{p}-1}),\forall
T>0$, because $u,J\in L^{\infty}(\mathcal{N}_{p,q,1}^{\frac{3}{p}-1})$. So,
given that we are still with the hypothesis $J_{0}=\nabla\times B_{0}$, it is
enough to repeat, but this time for all $T>0$, the same argument developed in
the proof of Theorem \ref{Result1} to conclude that $J=\nabla\times B$. Hence
it follows that $(u,B)$ is truly a global-in-time solution of system
(\ref{MHD}) in the sense of distributions.

To prove the uniqueness property, let us consider two solutions $(u^{1}%
,B^{1})$ and $(u^{2},B^{2})$ of system (\ref{MHD}), both starting from the
same initial data. We denote by $\Theta^{1}$ and $\Theta^{2}$ the
corresponding solutions of the extended system (\ref{System1})-(\ref{System2}%
). Also, without affecting generality, we can take $\Theta^{2}$ as the
previously constructed solution. Thus, it follows that%

\[
\Vert\Theta^{2}\Vert_{X}\leq c\Vert\Theta_{0}\Vert_{\mathcal{N}_{p,q,1}%
^{\frac{3}{p}-1}},\text{ \ \ if \ }\Vert\Theta_{0}\Vert_{\mathcal{N}%
_{p,q,1}^{\frac{3}{p}-1}}<\delta.
\]

Considering $\Phi=\Theta^{2}-\Theta^{1}$ and proceeding as in the proof of
Theorem \ref{Result1}, we have that $\Phi$ is given by (\ref{mildV}). Now,
making use of Lemma \ref{inhomogeneous-term}, we have, for all $T>0$, that%

\begin{align*}
&  \Vert\mathcal{B}(\Theta^{2},\Phi)\Vert_{L_{T}^{\infty}(\mathcal{N}%
_{p,q,1}^{\frac{3}{p}-1})}+\Vert\mathcal{B}(\Theta^{2},\Phi)\Vert_{L_{T}%
^{1}(\mathcal{N}_{p,q,1}^{\frac{3}{p}+1})}\\
&  \leq C\Vert\Pi(\Theta^{2},\Phi)\Vert_{L_{T}^{1}(\mathcal{N}_{p,q,1}%
^{\frac{3}{p}-1})}\\
&  =C\{\Vert\Pi_{a}(B^{2},\Phi_{2})-\Pi_{a}(u^{2},\Phi_{1})\Vert_{L_{T}%
^{1}(\mathcal{N}_{p,q,1}^{\frac{3}{p}-1})}+\Vert\Pi_{b}(B^{2},\Phi_{3}%
-\Phi_{1})\Vert_{L_{T}^{1}(\mathcal{N}_{p,q,1}^{\frac{3}{p}-1})}\\
&  +\Vert\nabla\times\Pi_{b}(curl^{-1}J^{2},\Phi_{3}-\Phi_{1})\Vert_{L_{T}%
^{1}(\mathcal{N}_{p,q,1}^{\frac{3}{p}-1})}\}.
\end{align*}
Hence, Lemma \ref{div(vxw)}, Lemma \ref{div(curlxw)} and Lemma
\ref{div(vxcurl)} yield%

\begin{align}
&  \Vert\mathcal{B}(\Theta^{2},\Phi)\Vert_{L_{T}^{\infty}(\mathcal{N}%
_{p,q,1}^{\frac{3}{p}-1})}+\Vert\mathcal{B}(\Theta^{2},\Phi)\Vert_{L_{T}%
^{1}(\mathcal{N}_{p,q,1}^{\frac{3}{p}+1})}\nonumber\\
&  \leq C\left\{  \int_{0}^{T}\Vert\Pi_{a}(B^{2},\Phi_{2})(s)-\Pi_{a}%
(u^{2},\Phi_{1})(s)\Vert_{\mathcal{N}_{p,q,1}^{\frac{3}{p}-1}}+\int_{0}%
^{T}\Vert\Pi_{b}(B^{2},\Phi_{3}-\Phi_{1})(s)\Vert_{\mathcal{N}_{p,q,1}%
^{\frac{3}{p}-1}}ds\right. \nonumber\\
&  \left.  +\int_{0}^{T}\Vert\Pi_{b}(curl^{-1}J^{2},\Phi_{3}-\Phi_{1}%
)(s)\Vert_{\mathcal{N}_{p,q,1}^{\frac{3}{p}}}ds\right\} \nonumber\\
&  \leq C\left\{  \int_{0}^{T}(\Vert B^{2}(s)\Vert_{\mathcal{N}_{p,q,1}%
^{\frac{3}{p}}}\Vert\Phi_{2}(s)\Vert_{\mathcal{N}_{p,q,1}^{\frac{3}{p}}}+\Vert
u^{2}(s)\Vert_{\mathcal{N}_{p,q,1}^{\frac{3}{p}}}\Vert\Phi_{1}(s)\Vert
_{\mathcal{N}_{p,q,1}^{\frac{3}{p}}})ds\right. \nonumber\\
&  \left.  +\int_{0}^{T}\Vert B^{2}(s)\Vert_{\mathcal{N}_{p,q,1}^{\frac{3}{p}%
}}\Vert(\Phi_{3}-\Phi_{1})(s)\Vert_{\mathcal{N}_{p,q,1}^{\frac{3}{p}}}%
ds+\int_{0}^{T}\Vert\Pi_{b}(curl^{-1}J^{2},\Phi_{3}-\Phi_{1})(s)\Vert
_{\mathcal{N}_{p,q,1}^{\frac{3}{p}}}ds\right\} \nonumber\\
&  \leq C\left\{  \int_{0}^{T}\Vert\Theta^{2}(s)\Vert_{\mathcal{N}%
_{p,q,1}^{\frac{3}{p}}}\Vert\Phi(s)\Vert_{\mathcal{N}_{p,q,1}^{\frac{3}{p}}%
}ds+\int_{0}^{T}\Vert\dv((curl^{-1}J^{2})\otimes(\Phi_{3}-\Phi_{1}%
))(s)\Vert_{\mathcal{N}_{p,q,1}^{\frac{3}{p}}}ds\right. \nonumber\\
&  \left.  +\int_{0}^{T}\Vert\dv((\Phi_{3}-\Phi_{1})\otimes(curl^{-1}%
J^{2}))(s)\Vert_{\mathcal{N}_{p,q,1}^{\frac{3}{p}}}ds\right\} \nonumber\\
&  \leq C\left\{  \int_{0}^{T}\Vert\Theta^{2}(s)\Vert_{\mathcal{N}%
_{p,q,1}^{\frac{3}{p}}}\Vert\Phi(s)\Vert_{\mathcal{N}_{p,q,1}^{\frac{3}{p}}%
}ds+\int_{0}^{T}\Vert J^{2}(s)\Vert_{\mathcal{N}_{p,q,1}^{\frac{3}{p}}}%
\Vert(\Phi_{3}-\Phi_{1})(s)\Vert_{\mathcal{N}_{p,q,1}^{\frac{3}{p}}}ds\right.
\nonumber\\
&  \left.  +\int_{0}^{T}\Vert J^{2}(s)\Vert_{\mathcal{N}_{p,q,1}^{\frac{3}%
{p}-1}}\Vert(\Phi_{3}-\Phi_{1})(s)\Vert_{\mathcal{N}_{p,q,1}^{\frac{3}{p}+1}%
}ds\right\} \nonumber\\
&  \leq C\left\{  \int_{0}^{T}\Vert\Theta^{2}(s)\Vert_{\mathcal{N}%
_{p,q,1}^{\frac{3}{p}}}\Vert\Phi(s)\Vert_{\mathcal{N}_{p,q,1}^{\frac{3}{p}}%
}ds+\int_{0}^{T}\Vert\Theta^{2}(s)\Vert_{\mathcal{N}_{p,q,1}^{\frac{3}{p}-1}%
}\Vert\Phi(s)\Vert_{\mathcal{N}_{p,q,1}^{\frac{3}{p}+1}}ds\right\}  .
\label{Uniestimate1}%
\end{align}

Applying an interpolation argument and Young inequality in (\ref{Uniestimate1}%
), we get, for all $\eta>0$, that%

\begin{align}
&  \Vert\mathcal{B}(\Theta^{2},\Phi)\Vert_{L_{T}^{\infty}(\mathcal{N}%
_{p,q,1}^{\frac{3}{p}-1})}+\Vert\mathcal{B}(\Theta^{2},\Phi)\Vert_{L_{T}%
^{1}(\mathcal{N}_{p,q,1}^{\frac{3}{p}+1})}\nonumber\\
&  \leq C\left\{  \int_{0}^{T}\Vert\Theta^{2}(s)\Vert_{\mathcal{N}%
_{p,q,1}^{\frac{3}{p}}}\Vert\Phi(s)\Vert_{\mathcal{N}_{p,q,1}^{\frac{3}{p}-1}%
}^{\frac{1}{2}}\Vert\Phi(s)\Vert_{\mathcal{N}_{p,q,1}^{\frac{3}{p}+1}}%
^{\frac{1}{2}}ds+\int_{0}^{T}\Vert\Theta^{2}(s)\Vert_{\mathcal{N}%
_{p,q,1}^{\frac{3}{p}-1}}\Vert\Phi(s)\Vert_{\mathcal{N}_{p,q,1}^{\frac{3}%
{p}+1}}ds\right\} \nonumber\\
&  \leq\eta\int_{0}^{T}\Vert\Phi(s)\Vert_{\mathcal{N}_{p,q,1}^{\frac{3}{p}+1}%
}ds+C\eta^{-1}\int_{0}^{T}\Vert\Theta^{2}(s)\Vert_{\mathcal{N}_{p,q,1}%
^{\frac{3}{p}}}^{2}\Vert\Phi(s)\Vert_{\mathcal{N}_{p,q,1}^{\frac{3}{p}-1}%
}ds+C\int_{0}^{T}\Vert\Theta^{2}(s)\Vert_{\mathcal{N}_{p,q,1}^{\frac{3}{p}-1}%
}\Vert\Phi(s)\Vert_{\mathcal{N}_{p,q,1}^{\frac{3}{p}+1}}ds\nonumber\\
&  \leq(\eta+C\Vert\Theta^{2}\Vert_{L_{T}^{\infty}(\mathcal{N}_{p,q,1}%
^{\frac{3}{p}-1})})\Vert\Phi\Vert_{L_{T}^{1}(\mathcal{N}_{p,q,1}^{\frac{3}%
{p}+1})}+C\eta^{-1}\int_{0}^{T}\Vert\Theta^{2}(s)\Vert_{\mathcal{N}%
_{p,q,1}^{\frac{3}{p}}}^{2}\Vert\Phi(s)\Vert_{\mathcal{N}_{p,q,1}^{\frac{3}%
{p}-1}}ds. \label{Uniestimate2}%
\end{align}
Using the same method as in obtaining (\ref{Uniestimate1}), we reach%

\begin{align*}
&  \Vert\mathcal{B}(\Phi,\Theta^{1})\Vert_{L_{T}^{\infty}(\mathcal{N}%
_{p,q,1}^{\frac{3}{p}-1})}+\Vert\mathcal{B}(\Phi,\Theta^{1})\Vert_{L_{T}%
^{1}(\mathcal{N}_{p,q,1}^{\frac{3}{p}+1})}\\
&  \leq C\left\{  \int_{0}^{T}\Vert\Theta^{1}(s)\Vert_{\mathcal{N}%
_{p,q,1}^{\frac{3}{p}}}\Vert\Phi(s)\Vert_{\mathcal{N}_{p,q,1}^{\frac{3}{p}}%
}ds+\int_{0}^{T}\Vert\Theta^{1}(s)\Vert_{\mathcal{N}_{p,q,1}^{\frac{3}{p}+1}%
}\Vert\Phi(s)\Vert_{\mathcal{N}_{p,q,1}^{\frac{3}{p}-1}}ds\right\}  .
\end{align*}
By employing the same procedure as utilized for (\ref{Uniestimate2}), we
derive, for all $\eta>0$, that%

\begin{align}
&  \Vert\mathcal{B}(\Phi,\Theta^{1})\Vert_{L_{T}^{\infty}(\mathcal{N}%
_{p,q,1}^{\frac{3}{p}-1})}+\Vert\mathcal{B}(\Phi,\Theta^{1})\Vert_{L_{T}%
^{1}(\mathcal{N}_{p,q,1}^{\frac{3}{p}+1})}\nonumber\\
&  \leq\eta\Vert\Phi\Vert_{L_{T}^{1}(\mathcal{N}_{p,q,1}^{\frac{3}{p}+1}%
)}+C\int_{0}^{T}(\Vert\Theta^{1}(s)\Vert_{\mathcal{N}_{p,q,1}^{\frac{3}{p}+1}%
}+\eta^{-1}\Vert\Theta^{1}(s)\Vert_{\mathcal{N}_{p,q,1}^{\frac{3}{p}}}%
^{2})\Vert\Phi(s)\Vert_{\mathcal{N}_{p,q,1}^{\frac{3}{p}-1}}ds.
\label{Uniestimate3}%
\end{align}
Now, it follows from equality (\ref{mildV}) and estimates (\ref{Uniestimate2})
and (\ref{Uniestimate3}) that%

\begin{align*}
&  \Vert\Phi\Vert_{L_{T}^{\infty}(\mathcal{N}_{p,q,1}^{\frac{3}{p}-1})}%
+\Vert\Phi\Vert_{L_{T}^{1}(\mathcal{N}_{p,q,1}^{\frac{3}{p}+1})}\\
&  \leq(2\eta+C\Vert\Theta^{2}\Vert_{L_{T}^{\infty}(\mathcal{N}_{p,q,1}%
^{\frac{3}{p}-1})})\Vert\Phi\Vert_{L_{T}^{1}(\mathcal{N}_{p,q,1}^{\frac{3}%
{p}+1})}\\
&  +C\eta^{-1}\int_{0}^{T}(\Vert\Theta^{1}(s)\Vert_{\mathcal{N}_{p,q,1}%
^{\frac{3}{p}+1}}+\Vert\Theta^{1}(s)\Vert_{\mathcal{N}_{p,q,1}^{\frac{3}{p}}%
}^{2}+\Vert\Theta^{2}(s)\Vert_{\mathcal{N}_{p,q,1}^{\frac{3}{p}}}^{2}%
)\Vert\Phi(s)\Vert_{\mathcal{N}_{p,q,1}^{\frac{3}{p}-1}}ds.
\end{align*}
And so, taking into account that
\[
\Vert\Theta^{2}\Vert_{L_{T}^{\infty}(\mathcal{N}_{p,q,1}^{\frac{3}{p}-1})}%
\leq\Vert\Theta^{2}\Vert_{L^{\infty}(\mathcal{N}_{p,q,1}^{\frac{3}{p}-1})}%
\leq\Vert\Theta^{2}\Vert_{X}\leq c\Vert\Theta_{0}\Vert_{\mathcal{N}%
_{p,q,1}^{\frac{3}{p}-1}}\text{ \ and \ }\Vert\Theta_{0}\Vert_{\mathcal{N}%
_{p,q,1}^{\frac{3}{p}-1}}<\delta,
\]
we have, for $\delta$ and $\eta$ small enough, that%

\[
\Vert\Phi\Vert_{L_{T}^{\infty}(\mathcal{N}_{p,q,1}^{\frac{3}{p}-1})}+\frac
{1}{2}\Vert\Phi\Vert_{L_{T}^{1}(\mathcal{N}_{p,q,1}^{\frac{3}{p}+1})}\leq
C\int_{0}^{T}(\Vert\Theta^{1}(s)\Vert_{\mathcal{N}_{p,q,1}^{\frac{3}{p}+1}%
}+\Vert\Theta^{1}(s)\Vert_{\mathcal{N}_{p,q,1}^{\frac{3}{p}}}^{2}+\Vert
\Theta^{2}(s)\Vert_{\mathcal{N}_{p,q,1}^{\frac{3}{p}}}^{2})\Vert\Phi
(s)\Vert_{\mathcal{N}_{p,q,1}^{\frac{3}{p}-1}}ds,
\]
which leads us to the conclusion that%

\[
\Vert\Phi(T)\Vert_{\mathcal{N}_{p,q,1}^{\frac{3}{p}-1}}\leq\Vert\Phi
\Vert_{L_{T}^{\infty}(\mathcal{N}_{p,q,1}^{\frac{3}{p}-1})}\leq C\int_{0}%
^{T}(\Vert\Theta^{1}(s)\Vert_{\mathcal{N}_{p,q,1}^{\frac{3}{p}+1}}+\Vert
\Theta^{1}(s)\Vert_{\mathcal{N}_{p,q,1}^{\frac{3}{p}}}^{2}+\Vert\Theta
^{2}(s)\Vert_{\mathcal{N}_{p,q,1}^{\frac{3}{p}}}^{2})\Vert\Phi(s)\Vert
_{\mathcal{N}_{p,q,1}^{\frac{3}{p}-1}}ds.
\]
By employing Gr\"{o}nwall inequality in the preceding estimate, we obtain
$\Vert\Phi(T)\Vert_{\mathcal{N}_{p,q,1}^{\frac{3}{p}-1}}=0,\forall T>0$, which
allows us to conclude that $\Theta_{1}\equiv\Theta_{2}.$ Our proof is complete.

\begin{flushright}
\ding{110}
\end{flushright}

\noindent\textbf{Acknowledgments.} LCFF was supported by CNPq (grant:
312484/2023-2), Brazil. RPS was supported by UTFPR, Brazil.

%%%%%%%%%%%%%%%%%%%%%%%%%%%%%%%%%%%%%%%%%%%%%%%%%%%%%%%%%%%%%%%%%
%% REFERÊNCIAS
%%%%%%%%%%%%%%%%%%%%%%%%%%%%%%%%%%%%%%%%%%%%%%%%%%%%%%%%%%%%%%%%

\end{document}